\documentclass[a4paper]{amsart}
\usepackage{graphicx,amsmath,amsfonts,latexsym,amssymb,amsthm,mathrsfs, color,hyperref}
\usepackage[latin1]{inputenc}
\newtheorem{thm}{Theorem}[section]

\newtheorem{corollary}[thm]{Corollary}
\newtheorem{definition}[thm]{Definition}
\newtheorem{theorem}[thm]{Theorem}
\newtheorem{lemma}[thm]{Lemma}
\newtheorem{remark}[thm]{Remark}

\begin{document}
\title[The fractional porous medium equation on conic manifolds]
{The fractional porous medium equation on manifolds with conical singularities I}
\author{Nikolaos Roidos and Yuanzhen Shao}
\address{Department of Mathematics, University of Patras, 26504 Rio Patras, Greece}
\email{roidos@math.upatras.gr}
\address{Department of Mathematics, The University of Alabama, Box 870350, Tuscaloosa, AL 35487-0350, USA}
\email{yshao8@ua.edu}

\thanks{\normalsize{The first author was supported by Deutsche Forschungsgemeinschaft, grant SCHR 319/9-1}}
\date{\today}
\begin{abstract} 
This is the first of a series of two papers which studies the fractional porous medium equation on a Riemannian manifold with isolated conical singularities. In this article, we show $R$-sectoriality for the fractional powers of possibly non-invertible $R$-sectorial operators. Applications concern existence, uniqueness and maximal $L^{q}$-regularity results for solutions of the fractional porous medium equation on manifolds with conical singularities. Space asymptotic behavior of the solutions close to the singularities is provided and its relation to the local geometry is established. Our method extends the freezing-of-coefficients method to the case of non-local operators that are expressed as linear combinations of terms in the form of a product of a function and a fractional power of a local operator.
\end{abstract}
\maketitle

\section{Introduction}

Let $X_{0}$ be a complex Banach space and let $A:\mathcal{D}(A)\rightarrow X_{0}$ be a closed linear operator that is sectorial of angle $\theta\in[0,\pi)$, i.e. the resolvent $(A+\lambda)^{-1}\in\mathcal{L}(X_{0})$ is defined for $\lambda\in\Lambda_{\theta}=\{z\in\mathbb{C}\backslash\{0\}\, |\, |\mathrm{arg}(z)|\leq\theta\}$ and moreover satisfies $|\lambda|\|(A+\lambda)^{-1}\|_{\mathcal{L}(X_{0})}\leq K$, $\lambda\in\Lambda_{\theta}$, for certain $K\geq1$; the class of such operators is denoted by $\mathcal{S}(\theta)$. In this case, by the functional calculus of sectorial operators, see e.g. \cite[Proposition III.4.6.10]{Am} or \cite[Section 15 C]{KW1} or \cite[Lemma 2.3.3]{Tan} or Theorem \ref{fracpower} below, for each $\sigma\in(0,1)$ the fractional power $A^{\sigma}$ of $A$ is a well defined closed linear operator in $X_{0}$ which is moreover sectorial of angle $\pi-(\pi-\theta)\sigma$.

In this paper we show that if in addition $A+c_{0}$ is $R$-sectorial of angle $\theta$ for certain $c_{0}\geq0$, i.e. if for each $\lambda_{1},\dots,\lambda_{N}\in \Lambda_{\theta}$, $x_{1},\dots,x_{N}\in X_{0}$, $N\in\mathbb{N}$, we have that
$$
\|\sum_{k=1}^{N}\epsilon_{k}\lambda_{k}(A+c_{0}+\lambda_{k})^{-1}x_{k}\|_{L^{2}(0,1;X_0)} \leq C \|\sum_{k=1}^{N}\epsilon_{k}x_{k}\|_{L^{2}(0,1;X_0)},
$$
for certain $C\geq1$ and the sequence of the Rademacher functions $\{\epsilon_{k}\}_{k\in\mathbb{N}}$, then there exists some $c\geq0$ such that $A^{\sigma}+c$ is $R$-sectorial of angle $\pi-(\pi-\theta)\sigma$. Hence, by denoting the class of $R$-sectorial operators of angle $\theta$ by $\mathcal{R}(\theta)$, our first result is the following.
\begin{theorem}\label{Rsec}
Let $\theta\in [0,\pi)$ and $A\in\mathcal{S}(\theta)$ such that $A+c_{0}\in\mathcal{R}(\theta)$ for certain $c_{0}\geq0$. Then for each $\sigma\in(0,1)$ there exists a $c\geq0$ such that $A^\sigma+c\in\mathcal{R}(\pi-(\pi-\theta)\sigma)${\em ;} in particular, we can choose $c=0$ when $c_{0}=0$.
\end{theorem}

It is well known that $R$-sectoriality is deeply related to the solvability and the regularity theory of linear and quasilinear parabolic problems, see e.g. \cite{DHP}, \cite{KaW}, \cite{PS} and \cite{Weis2}. Hence, in particular if the underlying space $X_{0}$ is UMD and $A$ has maximal $L^{q}$-regularity, i.e. the first order abstract linear Cauchy problem is well possed in the $L^{q}$-setting, see Section 2 for details, then due to standard theory, Theorem \ref{Rsec} implies that $A^{\sigma}$ has maximal $L^{q}$-regularity as well.

Next, as an application we consider a non-local evolution equation on manifolds with conical singularities. More precisely, let $\mathcal{B}$ be a smooth compact $(n+1)$-dimensional manifold, $n\ge 1$, with closed (i.e. compact without boundary) possibly disconnected smooth boundary $\partial\mathcal{B}$ of dimension $n$. We endow $\mathcal{B}$ with a degenerate Riemannian metric $g$ which in local coordinates $(x,y)\in[0,1)\times\partial\mathcal{B}$ on a collar neighborhood of the boundary is of the form $g= dx^{2}+x^{2}h$, where $h$ is a Riemannian metric on the cross-section $\partial\mathcal{B}$. We call $\mathbb{B}=(\mathcal{B},g)$ {\em manifold with conical singularities} or {\em conic manifold}; the boundary $\{0\}\times \partial\mathcal{B}$ of $\mathcal{B}$ corresponds to the conical tips. The Laplacian induced by $g$ on $(0,1)\times\partial\mathcal{B}$ has the degenerate form
\begin{equation}\label{DDelta}
\Delta=x^{-2}\big((x\partial_{x})^{2}+(n-1)(x\partial_{x})+\Delta_{h}\big),
\end{equation}
where $\Delta_{h}$ is the Laplacian on $\partial\mathbb{B}=(\partial\mathcal{B},h)$. 

We regard $\Delta$ as a second order cone differential operator acting on weighted Mellin-Sobolev spaces $\mathcal{H}_{p}^{s,\gamma}(\mathbb{B})$, $p\in(1,\infty)$, $s,\gamma\in\mathbb{R}$, see Definition \ref{MellSob}. It is well known that when $\Delta$ is considered as an unbounded operator in $\mathcal{H}_{p}^{s,\gamma}(\mathbb{B})$, it admits several closed extensions; each of these extensions corresponds to a subspace of a finite dimensional space $\mathcal{E}_{\Delta,\gamma}$ that is determined explicitly by the metric $h$, see Section 3 for details. Moreover, if we denote by $\mathbb{C}_{\omega}$ the space of smooth functions on $\mathbb{B}$ that are locally constant close to the singularities, see Definition \ref{constfunt}, it is known that under appropriate choice of the weight $\gamma$ in terms of the dimension and the local geometry, the map 
$$
\Delta: \mathcal{H}_{p}^{s+2,\gamma+2}(\mathbb{B})\oplus\mathbb{C}_{\omega}\rightarrow \mathcal{H}_{p}^{s,\gamma}(\mathbb{B})
$$
defines a closed extension $\underline{\Delta}_{s}$ of $\Delta$ in $\mathcal{H}_{p}^{s,\gamma}(\mathbb{B})$ such that $c_{0}-\underline{\Delta}_{s}\in \mathcal{R}(\theta)$ for each $c_{0}>0$ and $\theta\in[0,\pi)$, see \cite[Theorem 4.2]{RS2} or \cite[Theorem 6.7]{SS1}. By studying the nature of the pole zero of the resolvent of the above realization we show the following. 
 
\begin{theorem}\label{ThFracLap}
Let $p\in(1,\infty)$, $s\geq0$ and
\begin{equation}\label{gammaweight}
\frac{n-3}{2} <\gamma<\min\Big\{-1+\mu_{1},\frac{n+1}{2}\Big\}, \quad \mu_{j}=\sqrt{\left(\frac{n-1}2\right)^2-\lambda_{j}}, \quad j\in\mathbb{N}_{0}=\mathbb{N}\cup\{0\},
\end{equation}
where $\dots<\lambda_{1}<\lambda_{0}=0$ stands for the spectrum of $\Delta_{h}$. Moreover, consider the closed extension $\underline{\Delta}_{s}: \mathcal{H}_{p}^{s+2,\gamma+2}(\mathbb{B})\oplus\mathbb{C}_{\omega}\rightarrow \mathcal{H}_{p}^{s,\gamma}(\mathbb{B})$ of $\Delta$ in $\mathcal{H}_{p}^{s,\gamma}(\mathbb{B})$, where $\mathbb{C}_{\omega}$ denotes the space of smooth functions on $\mathbb{B}$ that are locally constant close to the singularities, see Definition \ref{constfunt}. Then, for each $\theta\in[0,\pi)$ we have that $-\underline{\Delta}_{s}\in \mathcal{S}(\theta)$. Therefore, for each $\sigma\in(0,1)$, through the functional calculus of sectorial operators, see e.g. Theorem \ref{fracpower}, the fractional power $(-\underline{\Delta}_{s})^{\sigma}:\mathcal{D}((-\underline{\Delta}_{s})^{\sigma})\rightarrow \mathcal{H}_{p}^{s,\gamma}(\mathbb{B})$ is a well defined closed linear operator that belongs to $ \mathcal{S}(\pi-(\pi-\theta)\sigma)$. The fractional Laplacian domain, described in Theorem \ref{fracpower}, satisfies 
\begin{equation}\label{domfrac}
\bigcup_{\varepsilon>0}\mathcal{H}_{p}^{s+2\sigma+\varepsilon,\gamma+2\sigma+\varepsilon}(\mathbb{B})\oplus\mathbb{C}_{\omega}\hookrightarrow \mathcal{D}((-\underline{\Delta}_{s})^{\sigma})\hookrightarrow \bigcap_{\varepsilon>0} \mathcal{H}_{p}^{s+2\sigma-\varepsilon,\gamma+2\sigma-\varepsilon}(\mathbb{B})\oplus\mathbb{C}_{\omega},
\end{equation}
and if in addition
\begin{equation}\label{gammapoles}
\gamma+2\sigma-1\notin\cup_{j\in\mathbb{N}_{0}}\{ \pm\mu_{j}\},
\end{equation}
then
\begin{equation}\label{sharpdom}
\mathcal{D}((-\underline{\Delta}_{s})^{\sigma})=\mathcal{H}_{p}^{s+2\sigma,\gamma+2\sigma}(\mathbb{B})\oplus\mathbb{C}_{\omega}.
\end{equation}
Furthermore, there exists a $c>0$ such that $(-\underline{\Delta}_{s})^{\sigma}+c\in\mathcal{R}(\pi-(\pi-\theta)\sigma)$. 
\end{theorem}

Note that the sum in \eqref{sharpdom} is either direct or we have $\mathbb{C}_{\omega}\subset \mathcal{H}_{p}^{s+2\sigma,\gamma+2\sigma}(\mathbb{B})$, so that for simplicity we use $\oplus$ instead of $+$.

Based on the above closed extension of the Laplacian we consider the fractional porous medium equation
\begin{eqnarray}\label{FPME1}
u'(t)+(-\Delta)^{\sigma}u^{m}(t) &=& 0, \quad t\in(0,T),\\\label{FPME2}
u(0) &=& u_{0},
\end{eqnarray}
where $\sigma\in(0,1)$, $m>0$, $T>0$ and $u_{0}$ is some given initial data. Here $(-\Delta)^{\sigma}$ is defined as a fractional power of a non-invertible sectorial operator as before. The problem \eqref{FPME1}-\eqref{FPME2} is a fractional version of the usual porous medium equation which is obtained after replacing $(-\Delta)^{\sigma}$ with $-\Delta$; concerning the usual porous medium equation, in order to avoid citing a large amount of literature, we only refer to the monograph \cite{Va} and to the references therein. Due to the non-locality of the fractional Laplacian $(-\Delta)^{\sigma}$, \eqref{FPME1}-\eqref{FPME2} can model long range diffusive interaction. As a consequence, the diffusion process described by \eqref{FPME1}-\eqref{FPME2} has applications to various fields, as heat control, statistical mechanics etc., see e.g. \cite{AtCa} and \cite{Jar}. 

The fractional porous medium equation has already been studied in $\mathbb{R}^{n}$ and the associated mathematical theory has been developed in several directions and under many aspects, see e.g. \cite{AtCa}, \cite{DPQRV0}, \cite{DPQRV1}, \cite{DPQRV2}, \cite{GMP}, \cite{GMP2} and \cite{Va1}. Note that in the above situations the fractional Laplacian is defined either through its Fourier transform symbol or by the self-adjointness of $\Delta$, so that it is always a particular case of a fractional power of a sectorial operator, see e.g. \cite[Theorem III.4.6.7]{Am}. Under this remark, in the present paper we present a different and more general approach to the problem \eqref{FPME1}-\eqref{FPME2} based on the maximal $L^{q}$-regularity theory for sectorial operators; the method we follow provides classical solutions with arbitrary high regularity and can be applied to various domains.

In Theorem \ref{wperturb} we show that for each strictly positive function $w$ that satisfies certain regularity, there exists a $c>0$ such that the operator $w(-\underline{\Delta}_{0})^{\sigma}+c$ is $R$-sectorial. Our method extends the standard freezing-of-coefficients method to the situation of non-local operators that are expressed as linear combinations of terms of the form $wA^{\sigma}$, where $A$ is a local operator. The key point here is the observation that a commutator of the form $[w,A^{\sigma}]$ is indeed of lower order in fractional sense.

In addition, in Theorem \ref{wperturbsect} we show that higher regularity in $w$ implies $R$-sectoriality of $w(-\underline{\Delta}_{s})^{\sigma}+c$ in higher order Mellin-Sobolev spaces. Here we use the non-commutative operator valued functional calculus theory for sectorial operators and, in particular, a theorem of Haller-Dintelmann and Hieber concerning the $\mathcal{H}^{\infty}$-calculus property for products of possibly non-commuting operators. 

The above two results show maximal $L^{q}$-regularity for the linearization of \eqref{FPME1} and are used for showing short time existence, uniqueness and maximal $L^{q}$-regularity for solutions of \eqref{FPME1}-\eqref{FPME2} by an abstract result of Cl\'ement and Li. Smoothness in time through the maximal $L^{q}$-regularity property is also shown by a theorem of Pr\"uss and Simonett. By denoting with $(\cdot,\cdot)_{\eta,q}$, $\eta\in(0,1)$, $q\in(1,\infty)$, the real interpolation functor of exponent $\eta$ and parameter $q$, we show the following well-posedness result for the fractional porous medium equation. 

\begin{theorem}\label{ThemonFrPME}
Let $\sigma_{0}=\max\{0,\frac{1}{2}(\frac{n+3}{2}-\mu_{1})\}<\sigma<1$, where $\mu_{1}$ is defined in \eqref{gammaweight}. Choose $p,q\in(1,\infty)$ such that 
$$
q>\frac{\sigma}{\sigma+\sigma_{0}} \quad \text{and} \quad \frac{n+1}{p}+\frac{2\sigma}{q}<2\sigma,
$$ 
and let 
$$
\gamma>\frac{n+1}{2}+\frac{2\sigma}{q}-2\sigma
$$ 
satisfying \eqref{gammaweight} and \eqref{gammapoles}. Then, the space $(\mathcal{H}_{p}^{2\sigma,\gamma+2\sigma}(\mathbb{B})\oplus\mathbb{C}_{\omega},\mathcal{H}_{p}^{0,\gamma}(\mathbb{B}))_{\frac{1}{q},q}$ consists of continuous functions on $\mathbb{B}$ and for each strictly positive 
\begin{equation}\label{u01}
u_{0}\in (\mathcal{H}_{p}^{2\sigma,\gamma+2\sigma}(\mathbb{B})\oplus\mathbb{C}_{\omega},\mathcal{H}_{p}^{0,\gamma}(\mathbb{B}))_{\frac{1}{q},q} \hookleftarrow \bigcup_{\varepsilon>0} \mathcal{H}_{p}^{2\sigma-\frac{2\sigma}{q}+\varepsilon,\gamma+2\sigma-\frac{2\sigma}{q}+\varepsilon}(\mathbb{B})\oplus\mathbb{C}_{\omega},
\end{equation}
there exists a $T>0$ and a unique 
\begin{equation}\label{regofu}
u\in W^{1,q}(0,T;\mathcal{H}_{p}^{s,\gamma}(\mathbb{B}))\cap L^{q}(0,T;\mathcal{H}_{p}^{s+2\sigma,\gamma+2\sigma}(\mathbb{B})\oplus\mathbb{C}_{\omega})
\end{equation}
solving \eqref{FPME1}-\eqref{FPME2}, where $s=0$. The solution also satisfies
\begin{eqnarray}\nonumber
\lefteqn{u\in C^{\infty}((0,T);\mathcal{H}_{p}^{2\sigma,\gamma+2\sigma}(\mathbb{B})\oplus\mathbb{C}_{\omega})}\\\label{rego55}
&&\cap\bigcap_{\varepsilon>0}C([0,T]; \mathcal{H}_{p}^{s+2\sigma-\frac{2\sigma}{q}-\varepsilon,\gamma+2\sigma-\frac{2\sigma}{q}-\varepsilon}(\mathbb{B})\oplus\mathbb{C}_{\omega})\hookrightarrow C([0,T];C(\mathbb{B})),
\end{eqnarray}
with $s=0$. If in particular
\begin{equation}\label{u02}
 u_{0}\in \bigcup_{\varepsilon>0}\mathcal{H}_{p}^{\nu+2+\frac{n+1}{p}+\varepsilon,\max\{\gamma+2,\frac{n+3}{2}\}+\varepsilon}(\mathbb{B})\oplus\mathbb{C}_{\omega}
\end{equation}
for some $\nu\geq0$, then the above $u$ satisfies \eqref{regofu}-\eqref{rego55} with $s=\nu$. 
\end{theorem}

We point out that the solution $u$ given by Theorem \ref{ThemonFrPME} is classical due to \eqref{sharpdom}. Moreover, the restriction of the fractional exponent $\sigma$ is necessary in our Mellin-Sobolev space setup. This is due to the non-linearity and due to the weight appearing in such spaces that describes the behavior of the functions close to the conical tips, in particular due to the necessity that elements in the interpolation space act by multiplication as bounded maps on the underlying space. 

Theorem \ref{ThemonFrPME} provides us information concerning the asymptotic behavior of the solution $u$ close to the singularities. More precisely, by \eqref{rego55} we can decompose the solution as $u=u_{\mathcal{H}}+u_{\mathcal{\mathbb{C}}}$, where $u_{\mathcal{H}}\in C^{\infty}((0,T); \mathcal{H}_{p}^{\nu+2\sigma,\gamma+2\sigma}(\mathbb{B}))$ and $u_{\mathcal{\mathbb{C}}}\in C^{\infty}((0,T); \mathbb{C}_{\omega})$. Moreover, by standard embedding properties of Mellin-Sobolev spaces, see e.g. \cite[Lemma 3.2]{RS3}, $u_{\mathcal{H}}\in C^{\infty}((0,T);C(\mathbb{B}))$ and, in local coordinates $(x,y)\in [0,1)\times\partial\mathcal{B}$ on the collar part, we have 
$$
|u_{\mathcal{H}}(t,x,y)|\leq c_{p}x^{\gamma+2\sigma-\frac{n+1}{2}}\|u_{\mathcal{H}}(t)\|_{\mathcal{H}_{p}^{\nu+2\sigma,\gamma+2\sigma}(\mathbb{B})}, \quad t\in (0,T), 
$$
where the constant $c_{p}>0$ depends only on $p$ and $\mathbb{B}$. Therefore, the Mellin-Sobolev part of the solution decays to zero close to the conical tips with certain rate that is determined by the local geometry, the initial data, the dimension and the fractional exponent.

Concerning the situation of the usual porous medium equation, the problem has already been considered on spaces with non-trivial geometry; we briefly mention the following contributions. In \cite{RS3} it was shown existence, uniqueness and maximal $L^{q}$-regularity for the short time solutions, where in \cite{RS4} this result was improved to long time existence and smoothness. Moreover, concerning the case of singular manifolds in the sense of H. Amann \cite{Am0}, in \cite{Sh1} it was shown existence, uniqueness and maximal continuous regularity for the short time solutions and in \cite{Sh2} global existence of $L^{1}$-mild solutions; see also \cite{Sh3} and \cite{Sh4} for similar problems on such spaces. For the case of the hyperbolic space, or more generally for Riemannian manifolds with nonpositive sectional curvature, we refer to \cite{GaMu}, \cite{GMV}, \cite{GMP3} and \cite{Va2}. 

The paper is organized as follows: Section 2 contains abstract theory concerning the fractional powers of possibly non-invertible sectorial operator as well as the maximal $L^{q}$-regularity property for linear and quasilinear parabolic problems; a proof of Theorem \ref{Rsec} is also included. In Section 3 we recall some basic theory of the naturally appearing differential operators on conic manifolds and, in particular, of the cone Laplacian. Section 4 is dedicated to the study of the model cone Laplacian, i.e. the analogue of $\Delta$ on the infinite cone $([0,\infty)\times\partial\mathcal{B},dx^{2}+x^{2}h)$. In Section 5 we prove Theorem \ref{ThFracLap}, i.e. that a particular realization of the cone Laplacian is sectorial, and therefore we can define its fractional powers. In Section 6 the fractional porous medium equation on conic manifolds is studied through the theory of maximal $L^{q}$-regularity and Theorem \ref{ThemonFrPME} is proved therein. Some elementary lemmas are collected in the Appendix on Section 7.

\section{Sectorial operators, functional calculus and maximal $L^{q}$-regularity}

Let $X_{1}\overset{d}{\hookrightarrow}X_{0}$ be a continuously and densely injected complex Banach couple. 

\begin{definition}[Sectoriality]\label{secdef}
Let $\mathcal{P}(K,\theta)$, $K\geq1$, $\theta\in[0,\pi)$, be the class of all closed densely defined linear operators $A$ in $X_{0}$ such that 
$$
S_{\theta}=\{\lambda\in\mathbb{C}\,|\, |\arg(\lambda)|\leq\theta\}\cup\{0\}\subset\rho{(-A)} \quad \mbox{and} \quad (1+|\lambda|)\|(A+\lambda)^{-1}\|_{\mathcal{L}(X_{0})}\leq K, \quad \lambda\in S_{\theta}.
$$
The elements in $\mathcal{P}(\theta)=\cup_{K\geq1}\mathcal{P}(K,\theta)$ are called {\em invertible sectorial operators of angle $\theta$} and for each $A\in\mathcal{P}(\theta)$ the constant $\inf\{K\, |\, A\in \mathcal{P}(K,\theta)\}$ is called {\em the sectorial bound of $A$}. 

Furthermore, denote by $\mathcal{S}(K,\theta)$ the supclass of $\mathcal{P}(K,\theta)$ such that if $A\in\mathcal{S}(K,\theta)$ then
$$
S_{\theta}\backslash\{0\}\subset\rho{(-A)} \quad \mbox{and} \quad |\lambda|\|(A+\lambda)^{-1}\|_{\mathcal{L}(E)}\leq K, \quad \lambda\in S_{\theta}\backslash\{0\}.
$$
The elements in $\mathcal{S}(\theta)=\cup_{K\geq1}\mathcal{S}(K,\theta)$ are called {\em sectorial operators of angle $\theta$} and for each $A\in\mathcal{S}(\theta)\backslash\mathcal{P}(\theta)$ the constant $\inf\{K\, |\, A\in \mathcal{S}(K,\theta)\}$ is called {\em the sectorial bound of $A$}. 
\end{definition}

Recall that $\mathcal{P}(K,\theta)\subset \mathcal{P}(2K+1,\phi)$ for certain $\phi\in(\theta,\pi)$, see e.g. \cite[(III.4.6.4)-(III.4.6.5)]{Am}, and similarly for the class $\mathcal{S}(\theta)$. Hence, whenever $A\in \mathcal{P}(\theta)$ or $A\in \mathcal{S}(\theta)$, we can always assume that $\theta>0$. Moreover, for any $\rho\geq0$ and $\theta\in(0,\pi)$, let $\Gamma_{\rho,\theta}$ be the counterclockwise oriented path defined by
$$
\Gamma_{\rho,\theta}=\{re^{-i\theta}\in\mathbb{C}\,|\,r\geq\rho\}\cup\{\rho e^{i\phi}\in\mathbb{C}\,|\,\theta\leq\phi\leq2\pi-\theta\}\cup\{re^{+i\theta}\in\mathbb{C}\,|\,r\geq\rho\}.
$$
We simply denote $\Gamma_{0,\theta}$ by $\Gamma_{\theta}$ and $\Gamma_{\theta}^{\pm}=\{re^{\pm i \theta}\in\mathbb{C}\, |\, r\geq0\}$. Furthermore, denote by $\Omega^{\circ}$ the interior of a domain $\Omega\subset\mathbb{C}$ and let $S_{\theta}^{\circ,\pm}=\{\lambda\in S_{\theta}^{\circ}\, |\, \pm\arg(\lambda)\geq0\}$.

The holomorphic functional calculus for sectorial operators in the class $\mathcal{P}(\theta)$ is defined by the Dunford integral formula, see e.g. \cite[Theorem 1.7]{DHP}. A typical example are the complex powers; for $\mathrm{Re}(z)<0$ they are defined by
\begin{equation}\label{cp}
A^{z}=\frac{1}{2\pi i}\int_{\Gamma_{\rho,\theta}}(-\lambda)^{z}(A+\lambda)^{-1}d\lambda,
\end{equation}
where $\rho>0$ is sufficiently small. The family $\{A^{z}\}_{\mathrm{Re}(z)<0}$ together with $A^{0}=I$ is a strongly continuous analytic semigroup on $X_{0}$, see e.g. \cite[Theorem III.4.6.2 and Theorem III.4.6.5]{Am}. Moreover, each $A^{z}$, $\mathrm{Re}(z)<0$, is injective and the complex powers for positive real part $A^{-z}$ are defined by $A^{-z}=(A^{z})^{-1}$, see e.g. \cite[(III.4.6.12)]{Am}. By Cauchy's theorem we can deform the path in \eqref{cp} and define the imaginary powers $A^{it}$, $t\in\mathbb{R}\backslash\{0\}$, as the closure of the operator
$$
A^{it}=\frac{\sin(i\pi t)}{i\pi t}\int_{0}^{\infty}s^{it}(A+s)^{-2}Ads \quad \text{in}\quad \mathcal{D}(A),
$$
see e.g. \cite[(III.4.6.21)]{Am}. For the properties of the complex powers of sectorial operators we refer to \cite[Theorem III.4.6.5]{Am}. Concerning the imaginary powers, suppose that that there exist some $\delta,M>0$ such that 
\begin{equation}\label{bipdef1}
A^{it}\in \mathcal{L}(X_{0}) \quad \text{and} \quad \|A^{it}\|_{\mathcal{L}(X_{0})}\leq M \quad \text{when} \quad t\in(-\delta,\delta).
\end{equation}
Then, see e.g. \cite[Theorem III.4.7.1 and Corollary III.4.7.2]{Am}, we have $A^{it}\in \mathcal{L}(X_{0})$ for each $t\in\mathbb{R}$ and there exist some $\varphi,\widetilde{M}>0$ such that 
\begin{equation}\label{bipdef2}
\|A^{it}\|_{\mathcal{L}(X_{0})}\leq \widetilde{M}e^{\varphi|t|}, \quad t\in\mathbb{R}.
\end{equation}

\begin{definition}[Bounded imaginary powers] Let $A\in\mathcal{P}(0)$ in $X_{0}$ such that \eqref{bipdef1}-\eqref{bipdef2} are satisfied. In this case we say that {\em $A$ has bounded imaginary powers} and denote $A\in\mathcal{BIP}(\varphi)$.
\end{definition}

The following property, stronger than the boundedness of the imaginary powers, can also be investigated for operators in the class $\mathcal{P}(\theta)$.

\begin{definition}[Bounded $H^{\infty}$-calculus]
Let $\theta\in(0,\pi)$, $\phi\in[0,\theta)$, $A\in\mathcal{P}(\theta)$ and let $H_{0}^{\infty}(\phi)$ be the space of all bounded analytic functions $f:\mathbb{C}\backslash S_{\phi}\rightarrow \mathbb{C}$ satisfying 
$$
|f(\lambda)|\leq c \Big(\frac{|\lambda|}{1+|\lambda|^{2}}\Big)^{\eta}\quad \text{for any} \quad \lambda\in \mathbb{C}\backslash S_{\phi} \quad \text{and some $c,\eta>0$ depending on $f$.}
$$
Any $f\in H_{0}^{\infty}(\phi)$ defines an element $f(-A)\in \mathcal{L}(X_{0})$ by 
\begin{equation}\label{hgsta}
f(-A)=\frac{1}{2\pi i}\int_{\Gamma_{\theta}}f(\lambda)(A+\lambda)^{-1} d\lambda.
\end{equation}
We say that the operator $A$ {\em has bounded $H^{\infty}$-calculus of angle $\phi$}, and we denote by $A\in \mathcal{H}^{\infty}(\phi)$, if there exists some $C>0$ such that
$$
\|f(-A)\|_{\mathcal{L}(X_{0})}\leq C\sup_{\lambda\in\mathbb{C}\backslash S_{\phi}}|f(\lambda)| \quad \mbox{for any} \quad f\in H_{0}^{\infty}(\phi).
$$
\end{definition}

We continue with the definition and the properties of the fractional powers of a possibly non-invertible sectorial operator in the class $\mathcal{S}(\theta)$.

\begin{theorem}[Fractional powers]\label{fracpower}
Let $\sigma\in(0,1)$, $\theta\in(0,\pi)$, $\phi\in[0,\theta)$ and $A\in\mathcal{S}(\theta)$ in $X_{0}$. Moreover, let
\begin{equation}\label{Resolvent1}
I_{\sigma}^{\pm}(\lambda)=\frac{\sin(\pi\sigma)}{\sigma}\int_{\Gamma_{\theta}^{\pm}}\frac{s^{\sigma}}{(s^{\sigma}+\lambda e^{i\pi\sigma})(s^{\sigma}+\lambda e^{-i\pi\sigma})}(A+s)^{-1}ds, \quad \lambda\in S_{\pi-(\pi-\phi)\sigma}^{\circ,\pm}.
\end{equation}
Then, there exists a unique $\sigma$-dependent closed linear operator $A^{\sigma}$ in $X_{0}$, called {\em $\sigma$-power of $A$}, such that $A^{\sigma}\in\mathcal{S}(\pi-(\pi-\phi)\sigma)$ and $(A^{\sigma}+\lambda)^{-1}=I_{\sigma}^{\pm}(\lambda)$ for all $\lambda\in S_{\pi-(\pi-\phi)\sigma}^{\circ,\pm}$; if $\lambda\in S_{\pi(1-\sigma)}^{\circ}$, then $I_{\sigma}^{\pm}(\lambda)=I_{\sigma}(\lambda)$, where
\begin{equation}\label{Resolvent2}
I_{\sigma}(\lambda)=\frac{\sin(\pi\sigma)}{\sigma}\int_{0}^{\infty}\frac{s^{\sigma}}{(s^{\sigma}+\lambda e^{i\pi\sigma})(s^{\sigma}+\lambda e^{-i\pi\sigma})}(A+s)^{-1}ds, \quad \lambda\in S_{\pi(1-\sigma)}^{\circ}.
\end{equation}
In particular, if $c>0$ then $(A+c)^{\sigma}$ is given by the usual Dunford integral formula
\begin{equation}\label{Aexp}
(A+c)^{\sigma}=\frac{\sin(\pi \sigma)}{\pi}\int_{0}^{\infty}s^{\sigma-1}(A+c)(A+c+s)^{-1}ds \quad \text{in} \quad \mathcal{D}(A),
\end{equation}
and $\mathcal{D}((A+c)^{\sigma})=\mathrm{Ran}((A+c)^{-\sigma})$ with
\begin{equation}\label{invpower}
(A+c)^{-\sigma}=\frac{\sin(\pi \sigma)}{\pi}\int_{0}^{\infty}s^{-\sigma}(A+c+s)^{-1}ds \in \mathcal{L}(X_{0}).
\end{equation}
Furthermore, 
\begin{equation}\label{compare}
\mathcal{D}((A+c)^{\sigma})=\mathcal{D}(A^{\sigma})\quad \text{and} \quad \|(A+c)^{\sigma}u-A^{\sigma}u\|_{X_{0}}\leq M c^{\sigma}\|u\|_{X_{0}}, \quad u\in \mathcal{D}(A^{\sigma}),
\end{equation}
for some $M$ depending only on $\sigma$ and the sectorial bound of $A$.
\end{theorem}
\begin{proof}
The above result is contained in \cite[Section 2.3.2]{Tan}; see \cite[Theorem 2.3.1 and Lermma 2.3.5]{Tan}. The integral formula representation \eqref{Resolvent1} for the resolvent can be seen as follows. Due to \cite[(2.40) and (2.44)]{Tan} we have that $(0,+\infty)\subset\rho(-A^{\sigma})$ and $(A^{\sigma}+s)^{-1}=I_{\sigma}(s)$ when $s\in(0,+\infty)$. Similarly to the proof of \cite[Proposition III.4.6.10]{Am}, for $\lambda,\lambda_{0}\in S_{\pi(1-\sigma)}^{\circ}$ we have
$$
I_{\sigma}(\lambda)-I_{\sigma}(\lambda_{0})=\frac{\sin(\pi\sigma)}{\sigma}\int_{0}^{\infty}\Theta(\lambda,\lambda_{0},s)(A+s)^{-1}ds,
$$
where 
$$
\Theta(\lambda,\lambda_{0},s)=\frac{(\lambda_{0}-\lambda)s^{\sigma}(2s^{\sigma}\cos(\pi\sigma)+\lambda+\lambda_{0})}{(s^{\sigma}+\lambda e^{i\pi\sigma})(s^{\sigma}+\lambda e^{-i\pi\sigma})(s^{\sigma}+\lambda_{0} e^{i\pi\sigma})(s^{\sigma}+\lambda_{0} e^{-i\pi\sigma})}.
$$
This shows the analyticity of $I_{\sigma}(\cdot)$ in $S_{\pi(1-\sigma)}^{\circ}$. 

{\em Extension argument}. For each $\delta\in(0,\pi(1-\sigma))$ there exists a $K_{\delta}>0$ such that 
$$
\|I_{\sigma}(\lambda)\|_{\mathcal{L}(X_{0})}\leq K_{\delta},\quad \lambda\in S_{\pi(1-\sigma)-\delta}\backslash\{\lambda\in\mathbb{C}\, |\, |\lambda|<\delta\}.
$$
Therefore, if we choose $r_{\delta}=(2K_{\delta})^{-1}$, by
\begin{equation}\label{steppf43}
(A^{\sigma}+\lambda)I_{\sigma}(\lambda_{0})=I+(\lambda-\lambda_{0})I_{\sigma}(\lambda_{0}),
\end{equation}
where $\lambda_{0}\in [\delta,\infty)$ and $|\lambda-\lambda_{0}|\leq r_{\delta}$, we deduce that $(A^{\sigma}+\lambda)^{-1}$ exists for each $\lambda$ in 
$$
 \Omega_{\delta}=\bigcup_{\lambda_{0}\in[\delta,\infty)}\{\lambda\in\mathbb{C}\, |\, |\lambda-\lambda_{0}|\leq r_{\delta}\}.
$$
The analyticity of $I_{\sigma}(\cdot)$ implies that $(A^{\sigma}+\lambda)^{-1}=I_{\sigma}(\lambda)$ for each $\lambda\in \Omega_{\delta}$. Hence, \eqref{steppf43} holds true even if $\lambda_{0}\in \Omega_{\delta}$ and $|\lambda-\lambda_{0}|\leq r_{\delta}$. After finitely many steps we can show that for each $\lambda_{0}\in S_{\pi(1-\sigma)-\delta}\backslash\{\lambda\in\mathbb{C}\, |\, |\lambda|<\delta\}$ we have $\{\lambda\in\mathbb{C}\, |\, |\lambda-\lambda_{0}|\leq r_{\delta}\}\subset\rho(-A^{\sigma})$ and $(A^{\sigma}+\mu)^{-1}=I_{\sigma}(\mu)$ when $\mu\in \{\lambda\in\mathbb{C}\, |\, |\lambda-\lambda_{0}|\leq r_{\delta}\}$. Due to the arbitrariness of $\delta$, we conclude that $S_{\pi(1-\sigma)}^{\circ}\subset\rho(-A^{\sigma})$ and $(A^{\sigma}+\lambda)^{-1}=I_{\sigma}(\lambda)$ when $\lambda\in S_{\pi(1-\sigma)}^{\circ}$.

Similarly for $\lambda, \lambda_{0}\in S_{\pi-(\pi-\phi)\sigma}^{\circ,\pm}$ we have
$$
I_{\sigma}^{\pm}(\lambda)-I_{\sigma}^{\pm}(\lambda_{0})=\frac{\sin(\pi\sigma)}{\sigma}\int_{\Gamma_{\theta}^{\pm}}\Theta(\lambda,\lambda_{0},s)(A+s)^{-1}ds,
$$
so that $I_{\sigma}^{\pm}(\cdot)$ is analytic in $S_{\pi-(\pi-\phi)\sigma}^{\circ,\pm}$. Moreover, if $\lambda\in S_{\pi(1-\sigma)}^{\circ}$, then we can deform the path of integration in $I_{\sigma}^{\pm}(\lambda)$ from $\Gamma_{\theta}^{\pm}$ to $[0,+\infty)$, so that $I_{\sigma}^{\pm}(\lambda)=I_{\sigma}(\lambda)$ when $\lambda\in S_{\pi(1-\sigma)}^{\circ}$. Finally, the extension argument above, applied with $S_{\pi(1-\sigma)-\delta}^{\circ}$ replaced by $S_{\pi-(\pi-\phi)\sigma-\delta}^{\circ,\pm}$ and $I_{\sigma}(\lambda)$ replaced by $I_{\sigma}^{\pm}(\lambda)$, shows that $S_{\pi-(\pi-\phi)\sigma}^{\circ,\pm}\subset\rho(-A^{\sigma})$ and $(A^{\sigma}+\lambda)^{-1}=I_{\sigma}^{\pm}(\lambda)$ when $\lambda\in S_{\pi-(\pi-\phi)\sigma}^{\circ,\pm}$.
\end{proof}

Consider the following abstract parabolic first order Cauchy problem
\begin{eqnarray}\label{app1}
u'(t)+Au(t)&=&f(t), \quad t\in(0,T),\\\label{app2}
u(0)&=&0,
\end{eqnarray}
where $-A:X_{1}\rightarrow X_{0}$ is the infinitesimal generator of an analytic semigroup on $X_{0}$ and $f\in L^q(0,T;X_{0})$, $q\in(1,\infty)$, $T>0$. The operator $A$ has {\em maximal $L^q$-regularity} if for any $f\in L^q(0,T;X_{0})$ there exists a unique $u\in W^{1,q}(0,T;X_{0})\cap L^{q}(0,T;X_{1})$ solving \eqref{app1}-\eqref{app2}; in this situation $u$ depends continuously on $f$ and the above property is independent of $q$ and $T$. 

\begin{definition}[$R$-boundedness]\label{rsec}
A set $E\subset \mathcal{L}(X_{0})$ is called {\em $R$-bounded} if for every $T_{1},\dots,T_{N}\in E$ and $x_{1},\dots,x_{N}\in X_0$, $N\in\mathbb{N}$, we have
$$
\|\sum_{k=1}^{N}\epsilon_{k}T_{k}x_{k}\|_{L^{2}(0,1;X_0)} \leq C \|\sum_{k=1}^{N}\epsilon_{k}x_{k}\|_{L^{2}(0,1;X_0)},
$$
for certain $C>0$, where $\{\epsilon_{k}\}_{k\in\mathbb{N}}$ is the sequence of Rademacher functions. The infimum of all such constants $C>0$ is called {\em the $R$-bound of $E$}.
\end{definition} 

According to the above definition, we recall next the notion of $R$-sectoriality; a boundedness property of the resolvent of a sectorial operator that is related to the maximal $L^{q}$-regularity. 

\begin{definition}[$R$-sectoriality]\label{rsec}
Denote by $\mathcal{R}(\theta)$, $\theta\in[0,\pi)$, the class of all operators $A\in \mathcal{S}(\theta)$ in $X_{0}$ such that the set $E=\{\lambda(A+\lambda)^{-1}\, |\, \lambda\in S_{\theta}\backslash\{0\}\}$ is $R$-bounded. If $A\in \mathcal{R}(\theta)$ then $A$ is called {\em $R$-sectorial of angle $\theta$} and the $R$-bound of $E$ is called {\em the $R$-sectorial bound of $A$}. 
\end{definition} 

If we restrict to the class of UMD (unconditionality of martingale differences property, see e.g. \cite[Section III.4.4]{Am}) Banach spaces, then we have the following.

\begin{theorem}[{\rm Kalton and Weis, \cite[Theorem 6.5]{KaW} or \cite[Theorem 4.2]{Weis2}}]\label{KaWeTh}
If $X_{0}$ is UMD and $A\in\mathcal{R}(\theta)$ in $X_{0}$ with $\theta>\pi/2$, then $A$ has maximal $L^{q}$-regularity. 
\end{theorem}

If an operator is $R$-sectorial then this property is passed to its fractional powers as we can see from the following. 

\subsubsection*{Proof of Theorem 1.1} By extending the area of $R$-sectoriality, see e.g. \cite[Section 4.1]{DHP}, we can assume that $\theta>0$ and that there exists some $\phi\in(\theta,\pi)$ such that $A\in\mathcal{R}(\phi)$.

(i) Assume first that $c_{0}=0$. Let 
$$
\{\lambda_{1},\dots,\lambda_{N}\}=\{r_{1}e^{i\psi_{1}},\dots,r_{N}e^{i\psi_{N}}\}\in S_{\pi-(\pi-\theta)\sigma}\backslash\{0\}
$$ 
and $x_{1},\dots,x_{N}\in X_{0}$, $N\in\mathbb{N}$. Denote by $\{\epsilon_{k}\}_{k\in\mathbb{N}}$ the sequence of the Rademacher functions and let $\phi(\psi_{k})=\mathrm{sign}(\psi_{k})\phi$ with the convention that $\phi(\psi_{k})=\phi$ when $\psi_{k}=0$. If $R_{A,\phi}$ is the $R$-sectorial bound of $A\in \mathcal{R}(\phi)$, then by \eqref{Resolvent1} we estimate
\begin{eqnarray*}
\lefteqn{\frac{\sigma}{\sin(\pi\sigma)}\|\sum_{k=1}^{N}\epsilon_{k}\lambda_{k}(A^{\sigma}+\lambda_{k})^{-1}x_{k}\|_{L^{2}(0,1;X_{0})}}\\
&=&\|\sum_{k=1}^{N}\epsilon_{k}\int_{0}^{\infty}\frac{s^{-\sigma}r_{k}e^{i(\psi_{k}+(1-\sigma)\phi(\psi_{k}))}}{(1+s^{-\sigma}r_{k}e^{i(\psi_{k}+(\pi-\phi(\psi_{k}))\sigma)})(1+s^{-\sigma}r_{k}e^{i(\psi_{k}-(\pi+\phi(\psi_{k}))\sigma)})}\\
&&\times(A+se^{i\phi(\psi_{k})})^{-1}x_{k}ds\|_{L^{2}(0,1;X_{0})}\\
&=&\frac{1}{\sigma}\|\sum_{k=1}^{N}\epsilon_{k}\int_{0}^{\infty}\frac{e^{i(\psi_{k}-\sigma\phi(\psi_{k}))}\big(\frac{r_{k}}{y}\big)^{\frac{1}{\sigma}}e^{i\phi(\psi_{k})}}{(1+ye^{i(\psi_{k}+(\pi-\phi(\psi_{k}))\sigma)})(1+ye^{i(\psi_{k}-(\pi+\phi(\psi_{k}))\sigma)})}\\
&&\times(A+\big(\frac{r_{k}}{y}\big)^{\frac{1}{\sigma}}e^{i\phi(\psi_{k})})^{-1}x_{k}dy\|_{L^{2}(0,1;X_{0})}\\
&\leq&\frac{R_{A,\phi}}{\sigma}\int_{0}^{\infty}\|\sum_{k=1}^{N}\epsilon_{k}\frac{e^{i(\psi_{k}-\sigma\phi(\psi_{k}))}}{(1+ye^{i(\psi_{k}+(\pi-\phi(\psi_{k}))\sigma)})(1+ye^{i(\psi_{k}-(\pi+\phi(\psi_{k}))\sigma)})}x_{k}\|_{L^{2}(0,1;X_{0})}dy\\
&\leq&2\frac{R_{A,\phi}}{\sigma}\Big(\int_{0}^{\infty}\sup_{\psi\leq|\pi-(\pi-\theta)\sigma|}\frac{1}{|1+ye^{i(\psi+(\pi-\phi(\psi))\sigma)}||1+ye^{i(\psi-(\pi+\phi(\psi))\sigma)}|}dy\Big)\\
&&\times \|\sum_{k=1}^{N}\epsilon_{k}x_{k}\|_{L^{2}(0,1;X_{0})},
\end{eqnarray*}
where at the last step we have used Kahane's contraction principle, see e.g. \cite[Proposition 2.5]{KW1}.

(ii) Suppose now that $c_{0}>0$. Clearly $A+c\in\mathcal{R}(\phi)$ for all $c\geq c_{0}$ and the $R$-sectorial bound of $A+c$ is uniformly bounded in $c$, see e.g. \cite[Lemma 2.6]{RS3}. Moreover, by the estimate in (i), the $R$-sectorial bound of $(A+c)^{\sigma}\in \mathcal{R}(\pi-(\pi-\theta)\sigma)$ is uniformly bounded in $c\geq c_{0}$; in particular the sectorial bound of $(A+c)^{\sigma}\in\mathcal{S}(0)$ is also uniformly bounded in $c\geq c_{0}$. The same holds true for $(A+c)^{\sigma}+c^{\sigma+\varepsilon}\in \mathcal{R}(\pi-(\pi-\theta)\sigma)$, where $\varepsilon>0$ is fixed. By \eqref{compare} we have that
$$
\|((A+c)^{\sigma}-A^{\sigma})((A+c)^{\sigma}+c^{\sigma+\varepsilon})^{-1}\|_{\mathcal{L}(X_{0})}\leq M c^{\sigma}\frac{K_{c}}{c^{\sigma+\varepsilon}},
$$
where $K_{c}$ is the sectorial bound of $(A+c)^{\sigma}\in\mathcal{S}(0)$ and $M$ depends only on $\sigma$ and the $\mathcal{S}(\phi)$-sectorial bound of $A$. Therefore, by writing 
$$
A^{\sigma}+c^{\sigma+\varepsilon}=(A+c)^{\sigma}+c^{\sigma+\varepsilon}+A^{\sigma}-(A+c)^{\sigma}
$$
and taking $c$ sufficiently large, we obtain the result by $R$-sectoriality perturbation, see e.g. \cite[Proposition 4.4.2]{PS}. \mbox{\ } \hfill $\square$

In practice sometimes we are interested in a subclass of $H^{\infty}$-calculus operators which satisfy the following stronger condition.

\begin{definition}[$R$-bounded $H^{\infty}$-calculus]\label{rsec}
Denote by $\mathcal{RH}(\theta)$, $\theta\in[0,\pi)$, the class of all operators $A\in \mathcal{H}(\theta)$ in $X_{0}$ such that the set $\{f(A)\, |\, f\in H_{0}^{\infty}(\theta), \sup_{\lambda\in\mathbb{C}\backslash S_{\theta}}|f(\lambda)|\leq1\}$ is $R$-bounded. Any $A\in \mathcal{RH}(\theta)$ is said to have {\em$R$-bounded $H^{\infty}$-calculus of angle $\theta$.} 
\end{definition} 

We recall a bounded $H^{\infty}$-calculus perturbation result for operators in the class $\mathcal{RH}(\theta)$; this will be used later for $R$-sectoriality perturbation. The result is obtained from non-commutative operator valued functional calculus theory of sectorial operators.

\begin{theorem}[{\rm Haller-Dintelmann and Hieber, \cite[Theorem 3.2]{DiHi}}]\label{HaHiTh} Let $A\in\mathcal{H}(\theta_{A})$, $B\in\mathcal{RH}(\theta_{B})$, $\theta_{A}+\theta_{B}>\pi$, such that $(B+\mu)^{-1}\mathcal{D}(A)\subseteq \mathcal{D}(A)$ for some (and hence for all) $\mu\in S_{\theta_{B}}$ and
\begin{equation}\label{DaPratoGrisv}
\|[A,(B+\mu)^{-1}](A+\lambda)^{-1}\|_{\mathcal{L}(X_{0})}\leq\frac{C}{(1+|\lambda|^{1-\alpha})(1+|\mu|^{1+\beta})}, \quad \lambda\in S_{\theta_{A}}, \mu\in S_{\theta_{B}},
\end{equation} 
for some $\alpha,C\geq0$ and $\beta>0$ satisfying $\alpha+\beta<1$. Then, for each $\theta\in[0,\theta_{A}+\theta_{B}-\pi)$ there exists a $c>0$ such that $AB+c$ with domain $\{u\in\mathcal{D}(B)\, |\, Bu\in\mathcal{D}(A)\}$ belongs to $\mathcal{H}(\theta)$.
\end{theorem}

Next, we describe an abstract maximal $L^{q}$-regularity result for quasilinear parabolic equations. Let $q\in(1,\infty)$, $U$ be an open subset of $(X_{1},X_{0})_{\frac{1}{q},q}$, $A(\cdot): U\rightarrow \mathcal{L}(X_{1},X_{0})$ and $F(\cdot,\cdot): U\times [0,T_{0}]\rightarrow X_{0}$, for some $T_{0}>0$. Consider the problem
\begin{eqnarray}\label{aqpp1}
u'(t)+A(u(t))u(t)&=&F(u(t),t)+G(t),\quad t\in(0,T),\\\label{aqpp2}
u(0)&=&u_{0},
\end{eqnarray}
where $T\in(0,T_{0})$, $u_{0}\in U$ and $G\in L^{q}(0,T_{0};X_{0})$. A Banach fixed point argument based on maximal $L^q$-regularity property for the linearization $A(u_{0})$ and on appropriate Lipschitz continuity conditions, implies the following short time existence result.
\begin{theorem}[{\rm Cl\'ement and Li, \cite[Theorem 2.1]{CL}}]\label{ClementLi}
Assume that:\\
{\em (H1)} $A(\cdot)\in C^{1-}(U;\mathcal{L}(X_{1},X_{0}))$.\\
{\em (H2)} $F(\cdot,\cdot)\in C^{1-,1-}(U\times [0,T_{0}];X_{0})$.\\
{\em (H3)} $A(u_{0})$ has maximal $L^q$-regularity.\\
Then, there exists a $T\in(0,T_{0})$ and a unique $u\in W^{1,q}(0,T;X_{0})\cap L^{q}(0,T;X_{1})$ solving \eqref{aqpp1}-\eqref{aqpp2}. 
\end{theorem}

Finally, we recall the following embedding of the maximal $L^q$-regularity space, namely
\begin{equation}\label{interpemb}
W^{1,q}(0,T;X_{0})\cap L^{q}(0,T;X_{1})\hookrightarrow C([0,T];(X_{1},X_{0})_{\frac{1}{q},q}), \quad q\in(1,\infty), \, T>0,
\end{equation}
see e.g. \cite[Theorem III.4.10.2]{Am}. 

\section{The Laplacian on a conic manifold}

 We regard $\Delta$ as {\em a cone differential operator} or {\em a Fuchs type operator} and recall some basic facts and results from the related underlying pseudo-differential theory, which is called {\em cone calculus}, towards the direction of the study of nonlinear partial differential equations. For more details we refer to \cite{CSS1}, \cite{GM}, \cite{GKM}, \cite{Kra}, \cite{Le}, \cite{Ro1}, \cite{RS2}, \cite{RS3}, \cite{RS4}, \cite{RS1}, \cite{SS1}, \cite{SS2}, \cite{SS}, \cite{Schu} and \cite{Sei}.

An $\mu$-th order, $\mu\in\mathbb{N}_{0}$, differential operator $A$ with smooth coefficients in the interior $\mathbb{B}^{\circ}$ of $\mathbb{B}$ is called a cone differential operator of order $\mu\in\mathbb{N}_{0}$ if its restriction to the collar part $(0,1)\times\partial\mathcal{B}$ admits the form 
\begin{equation}\label{Aconeww}
A=x^{-\mu}\sum_{k=0}^{\mu}a_{k}(x)(-x\partial_{x})^{k}, \quad \mbox{where} \quad a_{k}\in C^{\infty}([0,1);\mathrm{Diff}^{\mu-k}(\partial\mathbb{B})).
\end{equation}
Such an operator is called {\em $\mathbb{B}$-elliptic} if, in addition to the usual pseudodifferential symbol, its {\em rescaled symbol} (see e.g. \cite[(2.3)]{CSS1} for definition) is also pointwise invertible; this is the case for the Laplacian $\Delta$.

Cone differential operators act naturally on scales of {\em Mellin-Sobolev} spaces. Let $\omega\in C^{\infty}(\mathbb{B})$ be a fixed cut-off function near the boundary, i.e. a smooth non-negative function on $\mathcal{B}$ with $\omega=1$ near $\{0\}\times\partial \mathcal{B}$ and $\omega=0$ on $\mathcal{B}\backslash([0,1)\times \partial \mathcal{B})$. Moreover, assume that in local coordinates $(x,y)\in [0,1)\times \partial\mathcal{B}$, $\omega$ depends only on $x$. Denote by $C_{c}^{\infty}$ the space of smooth compactly supported functions and by $H_{p}^{s}$, $p\in(1,\infty)$, $s\in\mathbb{R}$, the usual Bessel potential space. 

\begin{definition}[Mellin-Sobolev spaces]\label{MellSob}
For any $\gamma\in\mathbb{R}$ consider the map 
$$
M_{\gamma}: C_{c}^{\infty}(\mathbb{R}_{+}\times\mathbb{R}^{n})\rightarrow C_{c}^{\infty}(\mathbb{R}^{n+1}) \quad \mbox{defined by} \quad u(x,y)\mapsto e^{(\gamma-\frac{n+1}{2})x}u(e^{-x},y). 
$$
Furthermore, take a covering $\kappa_{j}:U_{j}\subseteq\partial\mathcal{B} \rightarrow\mathbb{R}^{n}$, $j\in\{1,\dots,N\}$, $N\in\mathbb{N}$, of $\partial\mathcal{B}$ by coordinate charts and let $\{\phi_{j}\}_{j\in\{1,\dots,N\}}$ be a subordinate partition of unity. For any $p\in(1,\infty)$ and $s\in\mathbb{R}$ let $\mathcal{H}^{s,\gamma}_p(\mathbb{B})$ be the space of all distributions $u$ on $\mathbb{B}^{\circ}$ such that 
$$
\|u\|_{\mathcal{H}^{s,\gamma}_p(\mathbb{B})}=\sum_{j=1}^{N}\|M_{\gamma}(1\otimes \kappa_{j})_{\ast}(\omega\phi_{j} u)\|_{H^{s}_p(\mathbb{R}^{n+1})}+\|(1-\omega)u\|_{H^{s}_p(\mathbb{B})}
$$
is defined and finite, where $\ast$ refers to the push-forward of distributions. The space $\mathcal{H}^{s,\gamma}_p(\mathbb{B})$, called {\em (weighted) Mellin-Sobolev space}, is independent of the choice of the cut-off function $\omega$, the covering $\{\kappa_{j}\}_{j\in\{1,\dots,N\}}$ and the partition $\{\phi_{j}\}_{j\in\{1,\dots,N\}}${\em ;} if $A$ is as in \eqref{Aconeww}, then it induces a bounded map
$$
A: \mathcal{H}^{s+\mu,\gamma+\mu}_p(\mathbb{B}) \rightarrow \mathcal{H}^{s,\gamma}_p(\mathbb{B}).
$$
Finally, if $s\in \mathbb{N}_{0}$, then equivalently, $\mathcal{H}^{s,\gamma}_p(\mathbb{B})$ is the space of all functions $u$ in $H^s_{p,loc}(\mathbb{B}^\circ)$ such that, near the boundary, we have
$$
x^{\frac{n+1}2-\gamma}(x\partial_x)^{k}\partial_y^{\alpha}(\omega(x) u(x,y)) \in L_{loc}^{p}\big([0,1)\times \partial \mathcal{B}, \sqrt{\mathrm{det}[h]}\frac{dx}xdy\big),\quad k+|\alpha|\le s.
$$
\end{definition}

Note that since the usual Bessel potential spaces are UMD, by \cite[Theorem III.4.5.2]{Am}, the Mellin-Sobolev spaces are also UMD.

Next we restrict to the case of the Lapacian $\Delta$ and regard it as an unbounded operator in $\mathcal{H}^{s,\gamma}_p(\mathbb{B})$, $p\in(1,\infty)$, $s,\gamma\in\mathbb{R}$, with domain $C_{c}^{\infty}(\mathbb{B}^{\circ})$. The domain of its minimal extension (i.e. its closure) $\underline{\Delta}_{\min,s}$ is given by 
\begin{equation}\label{dminlap}
\mathcal{D}(\underline{\Delta}_{\min,s})=\Big\{u\in \bigcap_{\varepsilon>0}\mathcal{H}^{s+2,\gamma+2-\varepsilon}_p(\mathbb{B}) \, |\, \Delta u\in \mathcal{H}^{s,\gamma}_p(\mathbb{B})\Big\};
\end{equation}
in particular
$$
\mathcal{H}_{p}^{s+2,\gamma+2}(\mathbb{B})\hookrightarrow\mathcal{D}(\underline{\Delta}_{\min,s})\hookrightarrow \bigcap_{\varepsilon>0}\mathcal{H}_{p}^{s+2,\gamma+2-\varepsilon}(\mathbb{B}).
$$
If in addition the {\em conormal symbol} of $\Delta$, i.e. the following family of differential operators
$$
\mathbb{C} \ni \lambda \mapsto \lambda^{2}-(n-1)\lambda + \Delta_{h} \in \mathcal{L}(H_{2}^{2}(\partial\mathbb{B}),H_{2}^{0}(\partial\mathbb{B})),
$$
 is invertible on the line $\{\lambda\in\mathbb{C}\,|\, \mathrm{Re}(\lambda)= \frac{n-3}{2}-\gamma\}$, then we have precisely $\mathcal{D}(\underline{\Delta}_{\min,s})=\mathcal{H}^{s+2,\gamma+2}_p(\mathbb{B})$, i.e.
$$
\mathcal{D}(\underline{\Delta}_{\min,s})=\mathcal{H}_{p}^{s+2,\gamma+2}(\mathbb{B}) \quad \text{iff} \quad \pm\mu_{j}\neq \gamma+1, \, j\in\mathbb{N}_{0}.
$$ 

The domain of the maximal extension $\underline{\Delta}_{\max,s}$ of $\Delta$, defined by $\mathcal{D}(\underline{\Delta}_{\max,s})=\{u\in\mathcal{H}^{s,\gamma}_p(\mathbb{B}) \, |\, \Delta u\in \mathcal{H}^{s,\gamma}_p(\mathbb{B})\}$, is expressed as
\begin{equation}\label{dmax1}
\mathcal{D}(\underline{\Delta}_{\max,s})=\mathcal{D}(\underline{\Delta}_{\min,s})\oplus\mathcal{E}_{\Delta,\gamma}.
\end{equation}
Here
\begin{equation}\label{domflatcone}
\mathcal{E}_{\Delta,\gamma}=\bigoplus_{q_{j}^{\pm}\in I_{\gamma}} \mathcal{E}_{\Delta,\gamma,q_{j}^{\pm}}, \quad q_{j}^{\pm}=\frac{n-1}{2}\pm\mu_{j}, \quad j\in\mathbb{N}_{0}, \quad I_{\gamma}=\Big(\frac{n-3}{2}-\gamma,\frac{n+1}{2}-\gamma\Big),
\end{equation}
and for each $q_{j}^{\pm}$, $\mathcal{E}_{\Delta,\gamma,q_{j}^{\pm}}$ is a finite dimensional space consisting of $C^{\infty}(\mathbb{B}^{\circ})$-functions that vanish on $\mathcal{B}\backslash([0,1)\times \partial\mathcal{B})$ and in local coordinates on $(0,1)\times\partial\mathcal{B}$ they are of the form $\omega(x)c(y)x^{-q_{j}^{\pm}}\log^{k}(x)$, where $c\in C^{\infty}(\partial\mathbb{B})$ and $k\in\{0,1\}$. Note that $q_{j}^{\pm}$ are precisely the poles of the inverse of the conormal symbol of $\Delta$ and for each $q_{j}^{\pm}$ the exponent $k$ runs up to the order of the pole.

Due to \eqref{dmax1}, there are several closed extensions of $\Delta$ in $\mathcal{H}^{s,\gamma}_p(\mathbb{B})$; each one corresponds to a subspace of $\mathcal{E}_{\Delta,\gamma}$. For an overview on the domain structure of a general $\mathbb{B}$-elliptic cone differential operator we refer to \cite[Section 3]{GKM} or alternatively to \cite[Sections 2.2 and 2.3]{SS}. 

\begin{definition}\label{constfunt}
Recall that $\partial\mathcal{B}=\cup_{j=1}^{k_{\mathcal{B}}}\partial\mathcal{B}_{j}$, for certain $k_{\mathcal{B}}\in\mathbb{N}$, where $\partial\mathcal{B}_{j}$ are closed, smooth and connected. Denote by $\mathbb{C}_{\omega}$ the space of all $C^{\infty}(\mathbb{B}^\circ)$-functions $c$ that vanish on $\mathcal{B}\backslash([0,1)\times\partial\mathcal{B})$ and on each component $[0,1)\times\partial\mathcal{B}_{j}$, $j\in\{1,\dots,k_{\mathbb{B}}\}$, they are of the form $c_{j}\omega$, where $c_{j}\in\mathbb{C}$, i.e. $\mathbb{C}_{\omega}$ consists of smooth functions that are locally constant close to the boundary. Endow $\mathbb{C}_{\omega}$ with the norm $\|\cdot\|_{\mathbb{C}_{\omega}}$ given by $c\mapsto \|c\|_{\mathbb{C}_{\omega}}=(\sum_{j=1}^{k_{\mathcal{B}}}|c_{j}|^{2})^{\frac{1}{2}}$. 
\end{definition}

We close this section by recalling a particular close extension of $\Delta$. Under certain choice of the weight $\gamma$, $\mathbb{C}_{\omega}$ becomes a subspace of $\mathcal{E}_{\Delta,\gamma}$ and the realization of the Laplacian with domain $\mathcal{H}_{p}^{s+2,\gamma+2}(\mathbb{B})\oplus\mathbb{C}_{\omega}$ satisfies the property of maximal $L^q$-regularity as we can see from the following. 

\begin{theorem}\label{SecLapOrigin}
Let $p\in(1,\infty)$, $s\geq0$ and $\gamma$ be as in \eqref{gammaweight}. Consider the closed extension $\underline{\Delta}_{s}$ of the Laplacian $\Delta$ in 
$$
X_{0}^{s}=\mathcal{H}_{p}^{s,\gamma}(\mathbb{B})
$$ 
with domain 
\begin{equation}\label{delta1}
\mathcal{D}(\underline{\Delta}_{s})=X_{1}^{s}=\mathcal{H}_{p}^{s+2,\gamma+2}(\mathbb{B})\oplus\mathbb{C}_{\omega}.
\end{equation}
Then, for each $c>0$ and $\theta\in[0,\pi)$, we have $c-\underline{\Delta}_{s}\in\mathcal{R}(\theta)$.
\end{theorem}
\begin{proof}
This is \cite[Theorem 4.2]{RS2} together with \cite[Theorem 4.4.5]{PS}. See also \cite[Theorem 6.7]{SS1}.
\end{proof}

\section{The model cone Laplacian} 

Let us consider the differential operator
\begin{equation}\label{deltamodel}
\Delta_{\wedge}=x^{-2}\big((x\partial_{x})^{2}+(n-1)(x\partial_{x})+\Delta_{h}\big)
\end{equation}
acting on smooth functions on the infinite half cylinder 
$$
\partial\mathbb{B}^{\wedge}=([0,\infty)\times\partial\mathcal{B},dx^{2}+x^{2}h).
$$
$\Delta_{\wedge}$ is called {\em the model cone Laplacian} and $\partial\mathbb{B}^{\wedge}$ {\em the model cone} of $\mathbb{B}$.

\begin{definition}
Let $\kappa_{j}:U_{j}\subseteq\partial\mathcal{B} \rightarrow\mathbb{R}^{n}$, $j\in\{1,\dots,N\}$, $N\in\mathbb{N}$, be a covering of $\partial\mathcal{B}$ by coordinate charts and let $\{\phi_{j}\}_{j\in\{1,\dots,N\}}$ be a subordinate partition of unity. For any $p\in(1,\infty)$ and $s\in\mathbb{R}$ let $\mathcal{H}^{s}_{p,cone}(\mathbb{R}\times\partial\mathbb{B})$ be the space of all functions $u$ such that for each $j\in\{1,\dots,N\}$ we have
$$
(x,y)\mapsto \phi_{j}(\kappa_{j}^{-1}(y/\mathsf{x}(x)))u(x,\kappa_{j}^{-1}(y/\mathsf{x}(x))) \in H_{p}^{s}(\mathbb{R}\times\mathbb{R}^{n}),
$$ 
where $\mathbb{R}\ni x\mapsto\mathsf{x}(x)\in\mathbb{R}$ is a fixed smooth function that is equal to $x$ in $[-1/2,1/2]$, outside $[-1/2,1/2]$ is nonzero and outside $[-1,1]$ is constant. Moreover, if $\gamma\in\mathbb{R}$ let $\mathcal{K}_{p}^{s,\gamma}(\partial\mathbb{B}^{\wedge})$ be the space of all functions $v$ such that
$$
\omega v\in \mathcal{H}^{s,\gamma}_p(\mathbb{B}) \quad \text{and}\quad (1-\omega)v\in \mathcal{H}^{s}_{p,cone}(\mathbb{R}\times\partial\mathbb{B}).
$$
\end{definition}

The operator $\Delta_{\wedge}$ acts naturally on scales of Sobolev spaces $\mathcal{K}_{p}^{s,\gamma}(\partial\mathbb{B}^{\wedge})$, i.e. 
$$
\Delta_{\wedge}\in\mathcal{L}(\mathcal{K}_{p}^{s+2,\gamma+2}(\partial\mathbb{B}^{\wedge}),\mathcal{K}_{p}^{s,\gamma}(\partial\mathbb{B}^{\wedge})), \quad p\in(1,\infty), \, s,\gamma\in\mathbb{R}.
$$

\begin{remark}
Let $p\in(1,\infty)$ and $s,\gamma\in\mathbb{R}$. The scalar product in $\mathcal{H}^{0}_{2,cone}(\mathbb{R}\times\partial\mathbb{B})$ and $\mathcal{K}_{2}^{0,0}(\partial\mathbb{B}^{\wedge})$ identifies respectively the dual space of $\mathcal{H}^{s}_{p,cone}(\mathbb{R}\times\partial\mathbb{B})$ and $\mathcal{K}_{p}^{s,\gamma}(\partial\mathbb{B}^{\wedge})$ with $\mathcal{H}^{-s}_{p',cone}(\mathbb{R}\times\partial\mathbb{B})$ and $\mathcal{K}_{p'}^{-s,-\gamma}(\partial\mathbb{B}^{\wedge})$, where $1/p+1/p'=1$.
\end{remark}

Next, we show the following interpolation result concerning the spaces $\mathcal{H}^{s}_{p,cone}(\mathbb{R}\times\partial\mathbb{B})$.

\begin{lemma}\label{intHcone}
Let $p,q\in(1,\infty)$, $s\in\mathbb{R}$, $\rho>0$ and $\theta\in(0,1)$. For any $\varepsilon>0$ we have
$$
\mathcal{H}^{s+\rho\theta+\varepsilon}_{p,cone}(\mathbb{R}\times\partial\mathbb{B})\hookrightarrow(\mathcal{H}^{s}_{p,cone}(\mathbb{R}\times\partial\mathbb{B}),\mathcal{H}^{s+\rho}_{p,cone}(\mathbb{R}\times\partial\mathbb{B}))_{\theta,q}\hookrightarrow \mathcal{H}^{s+\rho\theta-\varepsilon}_{p,cone}(\mathbb{R}\times\partial\mathbb{B}).
$$
\end{lemma}
\begin{proof}
Let $\kappa_{j}:U_{j}\subseteq\partial\mathcal{B} \rightarrow\mathbb{R}^{n}$, $j\in\{1,\dots,N\}$, $N\in\mathbb{N}$, be a covering of $\partial\mathcal{B}$ by coordinate charts and let $\{\phi_{j}\}_{j\in\{1,\dots,N\}}$ be a subordinate partition of unity. If $u\in (\mathcal{H}^{s}_{p,cone}(\mathbb{R}\times\partial\mathbb{B}),\mathcal{H}^{s+\rho}_{p,cone}(\mathbb{R}\times\partial\mathbb{B}))_{\theta,q}$, then by the definition of the real interpolation we have
\begin{eqnarray*}
\lefteqn{\|u\|_{(\mathcal{H}^{s}_{p,cone}(\mathbb{R}\times\partial\mathbb{B}),\mathcal{H}^{s+\rho}_{p,cone}(\mathbb{R}\times\partial\mathbb{B}))_{\theta,q}}}\\
&=&\Big\|t^{-\theta}\inf\Big\{\sum_{j=1}^{N}\| \phi_{j}(\kappa_{j}^{-1}(y/\mathsf{x}))v_{1}(x,\kappa_{j}^{-1}(y/\mathsf{x}))\|_{H_{p}^{s}(\mathbb{R}\times\mathbb{R}^{n})}\\
&&+t\sum_{j=1}^{N}\| \phi_{j}(\kappa_{j}^{-1}(y/\mathsf{x}))v_{2}(x,\kappa_{j}^{-1}(y/\mathsf{x}))\|_{H_{p}^{s+\rho}(\mathbb{R}\times\mathbb{R}^{n})}\, \\
&&|\, v_{1}+v_{2}=u, v_{1}\in \mathcal{H}^{s}_{p,cone}(\mathbb{R}\times\partial\mathbb{B}), v_{2}\in \mathcal{H}^{s+\rho}_{p,cone}(\mathbb{R}\times\partial\mathbb{B})\Big\}\Big\|_{L^{q}(0,\infty;t^{-1}dt)}.
\end{eqnarray*}
Therefore, we estimate
\begin{eqnarray*}
N\lefteqn{\|u\|_{(\mathcal{H}^{s}_{p,cone}(\mathbb{R}\times\partial\mathbb{B}),\mathcal{H}^{s+\rho}_{p,cone}(\mathbb{R}\times\partial\mathbb{B}))_{\theta,q}}}\\
&\geq&\sum_{j=1}^{N}\Big\|t^{-\theta}\inf\Big\{\| \phi_{j}(\kappa_{j}^{-1}(y/\mathsf{x}))v_{1}(x,\kappa_{j}^{-1}(y/\mathsf{x}))\|_{H_{p}^{s}(\mathbb{R}\times\mathbb{R}^{n})}\\
&&+t\| \phi_{j}(\kappa_{j}^{-1}(y/\mathsf{x}))v_{2}(x,\kappa_{j}^{-1}(y/\mathsf{x}))\|_{H_{p}^{s+\rho}(\mathbb{R}\times\mathbb{R}^{n})}\, \\
&&|\, \phi_{j}(\kappa_{j}^{-1}(y/\mathsf{x}))v_{1}(x,\kappa_{j}^{-1}(y/\mathsf{x}))+\phi_{j}(\kappa_{j}^{-1}(y/\mathsf{x}))v_{2}(x,\kappa_{j}^{-1}(y/\mathsf{x}))\\
&&=\phi_{j}(\kappa_{j}^{-1}(y/\mathsf{x}))u(x,\kappa_{j}^{-1}(y/\mathsf{x})),\\
&& \phi_{j}(\kappa_{j}^{-1}(y/\mathsf{x}))v_{1}(x,\kappa_{j}^{-1}(y/\mathsf{x}))\in H_{p}^{s}(\mathbb{R}\times\mathbb{R}^{n}), \\
&&\phi_{j}(\kappa_{j}^{-1}(y/\mathsf{x}))v_{2}(x,\kappa_{j}^{-1}(y/\mathsf{x}))\in H_{p}^{s+\rho}(\mathbb{R}\times\mathbb{R}^{n})\Big\}\Big\|_{L^{q}(0,\infty;t^{-1}dt)}\\
&=&\sum_{j=1}^{N}\|\phi_{j}(\kappa_{j}^{-1}(y/\mathsf{x}))u(x,\kappa_{j}^{-1}(y/\mathsf{x}))\|_{(H_{p}^{s}(\mathbb{R}\times\mathbb{R}^{n}),H_{p}^{s+\rho}(\mathbb{R}\times\mathbb{R}^{n}))_{\theta,q}}\\
&\geq&C\sum_{j=1}^{N}\|\phi_{j}(\kappa_{j}^{-1}(y/\mathsf{x}))u(x,\kappa_{j}^{-1}(y/\mathsf{x}))\|_{H_{p}^{s+\rho\theta-\varepsilon}(\mathbb{R}\times\mathbb{R}^{n})},
\end{eqnarray*}
for certain $C>0$, where we have used the analogous result in $\mathbb{R}^{n+1}$, see e.g. \cite[(I.2.5.2)]{Am} together with \cite[Chapter 1, Theorem 7.1]{LiMa}. This shows the second embedding.

Since $C_{c}^{\infty}(\mathbb{R}\times\partial\mathbb{B})$ is dense in $\mathcal{H}^{s}_{p,cone}(\mathbb{R}\times\partial\mathbb{B})$ and $\mathcal{H}^{s+\rho}_{p,cone}(\mathbb{R}\times\partial\mathbb{B})$, by applying \cite[Section 1.11.2 (3a)]{Trib} to the above result, we obtain 
\begin{eqnarray*}
\lefteqn{\mathcal{H}^{-s-\rho\theta+\varepsilon}_{p',cone}(\mathbb{R}\times\partial\mathbb{B})}\\
&&\hookrightarrow (\mathcal{H}^{-s}_{p',cone}(\mathbb{R}\times\partial\mathbb{B}),\mathcal{H}^{-s-\rho}_{p',cone}(\mathbb{R}\times\partial\mathbb{B}))_{\theta,q'}=(\mathcal{H}^{-s-\rho}_{p',cone}(\mathbb{R}\times\partial\mathbb{B}),\mathcal{H}^{-s}_{p',cone}(\mathbb{R}\times\partial\mathbb{B}))_{1-\theta,q'},
\end{eqnarray*}
where $1/p+1/p'=1$, $1/q+1/q'=1$ and for the last equality we have used \cite[(I.2.5.4)]{Am}. Then, the first embedding follows by notting that $-s-\rho\theta+\varepsilon=-s-\rho+\rho(1-\theta)+\varepsilon$. 
\end{proof}

Similarly, concerning the spaces $\mathcal{K}_{p}^{s,\gamma}(\partial\mathbb{B}^{\wedge})$ we have the following.

\begin{lemma}\label{intKspace}
Let $p,q\in(1,\infty)$, $s,\gamma\in\mathbb{R}$, $\rho>0$ and $\theta\in(0,1)$. For any $\varepsilon>0$ we have
$$
\mathcal{K}^{s+\rho\theta+\varepsilon,\gamma+\rho\theta+\varepsilon}_p(\partial\mathbb{B}^{\wedge})\hookrightarrow
(\mathcal{K}^{s,\gamma}_p(\partial\mathbb{B}^{\wedge}),\mathcal{K}^{s+\rho,\gamma+\rho}_p(\partial\mathbb{B}^{\wedge}))_{\theta,q}
\hookrightarrow \mathcal{K}^{s+\rho\theta-\varepsilon,\gamma+\rho\theta-\varepsilon}_p(\partial\mathbb{B}^{\wedge}).
$$
\end{lemma}
\begin{proof}
For the first embedding, if $u\in (\mathcal{K}^{s,\gamma}_p(\partial\mathbb{B}^{\wedge}),\mathcal{K}^{s+\rho,\gamma+\rho}_p(\partial\mathbb{B}^{\wedge}))_{\theta,q}$, then
\begin{eqnarray*}
\lefteqn{\|u\|_{(\mathcal{K}^{s,\gamma}_p(\partial\mathbb{B}^{\wedge}),\mathcal{K}^{s+\rho,\gamma+\rho}_p(\partial\mathbb{B}^{\wedge}))_{\theta,q}}}\\
&=&\Big\|t^{-\theta}\inf\Big\{\|\omega u_{1}\|_{\mathcal{H}^{s,\gamma}_p(\mathbb{B})}+\|(1-\omega) u_{1}\|_{\mathcal{H}^{s}_{p,cone}(\mathbb{R}\times\partial\mathbb{B})}+t\|\omega u_{2}\|_{\mathcal{H}^{s+\rho,\gamma+\rho}_p(\mathbb{B})}\\
&&+t\|(1-\omega) u_{2}\|_{\mathcal{H}^{s+\rho}_{p,cone}(\mathbb{R}\times\partial\mathbb{B})}\\
&&|\, u_{1}+ u_{2}= u, u_{1}\in\mathcal{K}^{s,\gamma}_p(\partial\mathbb{B}^{\wedge}), u_{2}\in \mathcal{K}^{s+\rho,\gamma+\rho}_p(\partial\mathbb{B}^{\wedge}) \Big\}\Big\|_{L^{q}(0,+\infty;t^{-1}dt)}.
\end{eqnarray*}
Therefore
\begin{eqnarray*}
\lefteqn{\|u\|_{(\mathcal{K}^{s,\gamma}_p(\partial\mathbb{B}^{\wedge}),\mathcal{K}^{s+\rho,\gamma+\rho}_p(\partial\mathbb{B}^{\wedge}))_{\theta,q}}}\\
&\geq&\Big\|t^{-\theta}\inf\Big\{\|\omega u_{1}\|_{\mathcal{H}^{s,\gamma}_p(\mathbb{B})}+t\|\omega u_{2}\|_{\mathcal{H}^{s+\rho,\gamma+\rho}_p(\mathbb{B})}\, \\
&&|\, \omega(u_{1}+u_{2})=\omega u, \omega u_{1}\in\mathcal{H}^{s,\gamma}_p(\mathbb{B}), \omega u_{2}\in \mathcal{H}^{s+2,\gamma+2}_p(\mathbb{B}) \Big\}\Big\|_{L^{q}(0,+\infty;t^{-1}dt)}\\
&&+\Big\|t^{-\theta}\inf\Big\{\|(1-\omega) u_{1}\|_{\mathcal{H}^{s}_{p,cone}(\mathbb{R}\times\partial\mathbb{B})}+t\|(1-\omega) u_{2}\|_{\mathcal{H}^{s+\rho}_{p,cone}(\mathbb{R}\times\partial\mathbb{B})}\,|\, (1-\omega)(u_{1}+u_{2}-u)=0, \\
&&(1-\omega)u_{1}\in\mathcal{H}^{s}_{p,cone}(\mathbb{R}\times\partial\mathbb{B}), (1-\omega)u_{2}\in \mathcal{H}^{s+\rho,\gamma+\rho}_{p,cone}(\mathbb{R}\times\partial\mathbb{B}) \Big\}\Big\|_{L^{q}(0,+\infty;t^{-1}dt)}\\
&=&\|\omega u\|_{(\mathcal{H}^{s,\gamma}_p(\mathbb{B}),\mathcal{H}^{s+\rho,\gamma+\rho}_p(\mathbb{B}))_{\theta,q}}+\|(1-\omega)u\|_{(\mathcal{H}^{s}_{p,cone}(\mathbb{R}\times\partial\mathbb{B}),\mathcal{H}^{s+\rho}_{p,cone}(\mathbb{R}\times\partial\mathbb{B}))_{\theta,q}}\\
&\geq& C(\|\omega u\|_{\mathcal{H}^{s+\rho\theta-\varepsilon,\gamma+\rho\theta-\varepsilon}_p(\mathbb{B})}+\|(1-\omega)u\|_{\mathcal{H}^{s+\rho\theta-\varepsilon}_{p,cone}(\mathbb{R}\times\partial\mathbb{B})}),
\end{eqnarray*}
for certain $C>0$, where we have used \cite[Lemma 3.5]{RS3} and Lemma \ref{intHcone}.

Recall that $C_{c}^{\infty}(\partial\mathbb{B}^{\wedge})$ is dense in $\mathcal{K}^{s,\gamma}_p(\partial\mathbb{B}^{\wedge})$ and $\mathcal{K}^{s+\rho,\gamma+\rho}_p(\partial\mathbb{B}^{\wedge})$. By applying \cite[Section 1.11.2 (3a)]{Trib} to the above embedding and using \cite[(I.2.5.4)]{Am}, we get
\begin{eqnarray*}
\lefteqn{\mathcal{K}^{-s-\rho\theta+\varepsilon,-\gamma-\rho\theta+\varepsilon}_{p'}(\partial\mathbb{B}^{\wedge})}\\
&&\hookrightarrow (\mathcal{K}^{-s,-\gamma}_{p'}(\partial\mathbb{B}^{\wedge}),\mathcal{K}^{-s-\rho,-\gamma-\rho}_{p'}(\partial\mathbb{B}^{\wedge}))_{\theta,q'}=(\mathcal{K}^{-s-\rho,-\gamma-\rho}_{p'}(\partial\mathbb{B}^{\wedge}),\mathcal{K}^{-s,-\gamma}_{p'}(\partial\mathbb{B}^{\wedge}))_{1-\theta,q'},
\end{eqnarray*}
where as usual $1/p+1/p'=1$ and $1/q+1/q'=1$. Then, the result follows since $-s-\rho\theta+\varepsilon=-s-\rho+\rho(1-\theta)+\varepsilon$ and $-\gamma-\rho\theta+\varepsilon=-\gamma-\rho+\rho(1-\theta)+\varepsilon$. 
\end{proof}

Let us now consider $\Delta_{\wedge}$ as an unbounded operator in $\mathcal{K}_{p}^{s,\gamma}(\partial\mathbb{B}^{\wedge})$, $p\in(1,\infty)$, $s,\gamma\in\mathbb{R}$, with domain $C_{c}^{\infty}(\partial\mathbb{B}^{\wedge})$. The domain of its maximal extension $\underline{\Delta}_{\wedge,\max,s}$ differs from the domain of its minimal extension $\underline{\Delta}_{\wedge,\min,s}$ by an $s$-independent finite dimensional space $\mathcal{E}_{\Delta,\gamma}^{\wedge}$, which is also called {\em asymptotics space}; $\mathcal{E}_{\Delta,\gamma}^{\wedge}$ is isomorphic to $\mathcal{E}_{\Delta,\gamma}$ in \eqref{dmax1} and has similar structure, see e.g. \cite[Proposition 2.11]{SS}, \cite[Section 3]{SS1} or \cite[Theorem 4.7]{GKM}. More precisely, we have
\begin{equation}\label{dmaxmodelcone}
\mathcal{D}(\underline{\Delta}_{\wedge,\max,s})=\mathcal{D}(\underline{\Delta}_{\wedge,\min,s})\oplus\mathcal{E}_{\Delta,\gamma}^{\wedge}=\mathcal{D}(\underline{\Delta}_{\wedge,\min,s})\oplus\bigoplus_{q_{j}^\pm\in I_{\gamma}}\mathcal{E}_{\Delta,\gamma,q_{j}^{\pm}}^{\wedge}.
\end{equation}
Here, for the domain of the closure we have
$$
\mathcal{D}(\underline{\Delta}_{\wedge,\min,s})=\Big\{u\in \bigcap_{\varepsilon>0}\mathcal{K}_{p}^{s+2,\gamma+2-\varepsilon}(\partial\mathbb{B}^{\wedge})\, |\, \Delta_{\wedge}u\in \mathcal{K}_{p}^{s,\gamma}(\partial\mathbb{B}^{\wedge})\Big\};
$$
in particular
$$
\mathcal{K}_{p}^{s+2,\gamma+2}(\partial\mathbb{B}^{\wedge})\hookrightarrow\mathcal{D}(\underline{\Delta}_{\wedge,\min,s})\hookrightarrow \bigcap_{\varepsilon>0}\mathcal{K}_{p}^{s+2,\gamma+2-\varepsilon}(\partial\mathbb{B}^{\wedge})
$$
and
\begin{equation}\label{dminmodcone}
\mathcal{D}(\underline{\Delta}_{\wedge,\min,s})=\mathcal{K}_{p}^{s+2,\gamma+2}(\partial\mathbb{B}^{\wedge}) \quad \text{iff} \quad \pm\mu_{j}\neq \gamma+1, \, j\in\mathbb{N}_{0}.
\end{equation}
Moreover, for each $q_{j}^{\pm}$, which is given by \eqref{domflatcone}, $\mathcal{E}_{\Delta,\gamma,q_{j}^{\pm}}^{\wedge}$ is a finite dimensional space consisting of $C^{\infty}((\partial\mathbb{B}^{\wedge})^{\circ})$-functions that vanish on $[1,\infty)\times \partial\mathcal{B}$ and in local coordinates on $(0,1)\times\partial\mathcal{B}$ they are of the form $\omega(x)c(y)x^{-q_{j}^{\pm}}\log^{k}(x)$, where $c\in C^{\infty}(\partial\mathbb{B})$ and $k\in\{0,1\}$.

\begin{lemma}\label{inter} 
Let $p,q\in(1,\infty)$, $s\in\mathbb{R}$, $\gamma\in (\frac{n-3}{2},\frac{n+1}{2})$ and $\theta\in(0,1)$. Then, the following embeddings hold
\begin{eqnarray*}
\lefteqn{\mathcal{K}^{s+2\theta+\varepsilon,\gamma+2\theta+\varepsilon}_p(\partial\mathbb{B}^{\wedge})\oplus\mathbb{C}_{\omega}}\\
&&\hookrightarrow(\mathcal{K}^{s,\gamma}_p(\partial\mathbb{B}^{\wedge}),\mathcal{K}^{s+2,\gamma+2}_p(\partial\mathbb{B}^{\wedge})\oplus\mathbb{C}_{\omega})_{\theta,q} \hookrightarrow \mathcal{K}^{s+2\theta-\varepsilon,\gamma+2\theta-\varepsilon}_p(\partial\mathbb{B}^{\wedge})\oplus\mathbb{C}_{\omega},
\end{eqnarray*}
for every $\varepsilon>0$. 
\end{lemma}
\begin{proof}
Concerning the first embedding, by standard properties of interpolation spaces (see e.g. \cite[Proposition I.2.3.2]{Am}) we have that
$$
(\mathcal{K}^{s,\gamma}_p(\partial\mathbb{B}^{\wedge}),\mathcal{K}^{s+2,\gamma+2}_p(\partial\mathbb{B}^{\wedge}))_{\theta,q}\hookrightarrow(\mathcal{K}^{s,\gamma}_p(\partial\mathbb{B}^{\wedge}),\mathcal{K}^{s+2,\gamma+2}_p(\partial\mathbb{B}^{\wedge})\oplus\mathbb{C}_{\omega})_{\theta,q}.
$$
Therefore, by Lemma \ref{intKspace} we obtain
$$
\mathcal{K}^{s+2\theta+\varepsilon,\gamma+2\theta+\varepsilon}_p(\partial\mathbb{B}^{\wedge})\hookrightarrow(\mathcal{K}^{s,\gamma}_p(\partial\mathbb{B}^{\wedge}),\mathcal{K}^{s+2,\gamma+2}_p(\partial\mathbb{B}^{\wedge})\oplus\mathbb{C}_{\omega})_{\theta,q},
$$
and the result follows by
$$
\mathbb{C}_{\omega}\hookrightarrow(\mathcal{K}^{s,\gamma}_p(\partial\mathbb{B}^{\wedge}),\mathcal{K}^{s+2,\gamma+2}_p(\partial\mathbb{B}^{\wedge})\oplus\mathbb{C}_{\omega})_{\theta,q}.
$$

Concerning the second embedding, if $u\in(\mathcal{K}^{s,\gamma}_p(\partial\mathbb{B}^{\wedge}),\mathcal{K}^{s+2,\gamma+2}_p(\partial\mathbb{B}^{\wedge})\oplus\mathbb{C}_{\omega})_{\theta,q}$, then in local coordinates $(x,y_{1},\dots,y_{n})\in [0,\infty)\times \partial\mathcal{B}$ we have that $(x\partial_x)^{2} u$, $x\partial_x u$ and $\partial_{y_{i}}\partial_{y_{j}}u$, $i,j\in\{1,..,n\}$, belong to 
\begin{equation}\label{conddx}
(\mathcal{K}^{s-2,\gamma}_p(\partial\mathbb{B}^{\wedge}),\mathcal{K}^{s,\gamma+2}_p(\partial\mathbb{B}^{\wedge}))_{\theta,q}\hookrightarrow \mathcal{K}^{s+2\theta-2-\varepsilon,\gamma+2\theta-\varepsilon}_p(\partial\mathbb{B}^{\wedge}),
\end{equation}
for any $\varepsilon>0$, where we have used Lemma \ref{intKspace}. Therefore, by \eqref{deltamodel} $u$ belongs to the maximal domain of $\Delta_{\wedge}$ in $\mathcal{K}^{s+2\theta-2-\varepsilon,\gamma+2\theta-2-\varepsilon}_p(\partial\mathbb{B}^{\wedge})$, i.e. to the right hand side of \eqref{dmaxmodelcone} with $s$ replaced by $s+2(\theta-1)-\varepsilon$ and $\gamma$ replaced by $\gamma+2(\theta-1)-\varepsilon$. Since the sum in \eqref{dmaxmodelcone} is direct, by the formula 
$$
x\partial_{x}(x^{-q_{j}^{\pm}}\log^{k}(x))=-q_{j}^{\pm}x^{-q_{j}^{\pm}}\log^{k}(x)+kx^{-q_{j}^{\pm}}\log^{k-1}(x), \quad k\in\{0,1\}, 
$$
together with the regularity of $(x\partial_x)^{2} u$ given by \eqref{conddx} we conclude that the only asymptotics space component contributing to the expression of $u$ is the one that corresponds to $q_j^\pm=0$, i.e. the space $\mathbb{C}_{\omega}$.
\end{proof}

Next, we focus on {\em dilation invariant extensions} of the model cone Laplacian, as e.g. in \cite[Section 3.2 (E2)]{SS}. It is well known that, if the spectrum is contained in a sector, then such an extension is necessarily sectorial.

\begin{lemma}\label{Rsecmodcone}
Let $p\in(1,\infty)$, $\gamma\in\mathbb{R}$ such that \eqref{dminmodcone} is satisfied and, according to \eqref{dmaxmodelcone}, consider a closed extension $\underline{\Delta}_{\wedge}$ of the model cone Laplacian $\Delta_{\wedge}$ in $\mathcal{K}_{p}^{0,\gamma}(\partial\mathbb{B}^{\wedge})$ with domain $\mathcal{D}(\underline{\Delta}_{\wedge})=\mathcal{K}_{p}^{2,\gamma+2}(\partial\mathbb{B}^{\wedge})\oplus \underline{\mathcal{E}}_{\Delta,\gamma}^{\wedge}$, where $\underline{\mathcal{E}}_{\Delta,\gamma}^{\wedge}$ is a subspace of $\mathcal{E}_{\Delta,\gamma}^{\wedge}$. Assume that:\\
{\em (i)} $\underline{\mathcal{E}}_{\Delta,\gamma}^{\wedge}$ is invariant under dilations, i.e. in local coordinates $(x,y)\in (0,\infty)\times\partial\mathcal{B}$ we have that, if $u(x,y)\in \underline{\mathcal{E}}_{\Delta,\gamma}^{\wedge}$, then $u(\rho x,y)\in \underline{\mathcal{E}}_{\Delta,\gamma}^{\wedge}$ for each $\rho>0$. \\
{\em (ii)}The spectrum of $\underline{\Delta}_{\wedge}$ is contained in $\mathbb{C}\backslash S_{\phi}^{\circ}$, for some $\phi\in (0,\pi)$.\\
Then, $-\underline{\Delta}_{\wedge}\in\mathcal{S}(\theta)$ for any $\theta\in[0,\phi)$. 
\end{lemma}
\begin{proof}
Let $\kappa_{\rho}$, $\rho>0$, be the normalized dilation group action on functions on $\partial\mathbb{B}^{\wedge}$ defined by $(\kappa_{\rho}u)(x,y)=\rho^{\eta} u(\rho x,y)$, $(x,y)\in(0,\infty)\times\partial\mathcal{B}$, where $\eta=\frac{n+1}{2}-\gamma$, see e.g. \cite[Definition 2.3]{GKM} or \cite[Definition 5.12]{Kra}. Then, similarly to \cite[(7.1)]{GKM} or \cite[(5.18)]{Kra}, we have that
\begin{equation}\label{modelcone}
\lambda-\Delta_{\wedge}=\rho^{2}\kappa_{\rho}(\rho^{-2}\lambda-\Delta_{\wedge})\kappa_{\rho}^{-1}, \quad \lambda\in \mathbb{C}, \, \rho>0.
\end{equation}
Note that $\kappa_{\rho}$ is an isometry on $\mathcal{K}_{p}^{0,\gamma}(\partial\mathbb{B}^{\wedge})$ and moreover $\mathcal{D}(\underline{\Delta}_{\wedge})$ is invariant under $\kappa_{\rho}$. Therefore, \eqref{modelcone} still holds if we replace $\Delta_{\wedge}$ with $\underline{\Delta}_{\wedge}$. Thus, by taking $\rho=\sqrt{|\lambda|}$ we obtain
$$
|\lambda|(\lambda-\underline{\Delta}_{\wedge})^{-1}=\kappa_{\rho}(|\lambda|^{-1}\lambda-\underline{\Delta}_{\wedge})^{-1}\kappa_{\rho}^{-1}, \quad \lambda\in S_{\theta},
$$
and hence $-\underline{\Delta}_{\wedge}\in\mathcal{S}(\theta)$. 
\end{proof}

In particular, the model cone analogue of the closed extension \eqref{delta1} is sectorial.

\begin{corollary}\label{SecLap}
Let $p\in(1,\infty)$ and $\gamma$ be as in \eqref{gammaweight}. According to \eqref{dmaxmodelcone}, consider the closed extension $\underline{\Delta}_{\wedge}$ of the model cone Laplacian $\Delta_{\wedge}$ in $\mathcal{K}_{p}^{0,\gamma}(\partial\mathbb{B}^{\wedge})$ with domain 
\begin{equation}\label{delta}
\mathcal{D}(\underline{\Delta}_{\wedge})=\mathcal{K}_{p}^{2,\gamma+2}(\partial\mathbb{B}^{\wedge})\oplus \mathbb{C}_{\omega}.
\end{equation}
Then, $-\underline{\Delta}_{\wedge}\in\mathcal{S}(\theta)$ for any $\theta\in[0,\pi)$. In particular, $0$ is a simple pole of $(\lambda-\underline{\Delta}_{\wedge})^{-1}$.
\end{corollary}
\begin{proof}
By returning to the conic manifold $\mathbb{B}$, the closed extension $\underline{\Delta}_{0}$ from Theorem \ref{SecLapOrigin} satisfies the assumptions of \cite[Theorem 5.6]{SS} and the conditions (i), (ii) and (iii) of \cite[Theorem 5.7]{SS}. Therefore, according to \cite[Theorem 5.6 and Theorem 5.7]{SS}, see also \cite[Theorem 2.9 and Remark 2.10]{RS1}, the closed extension $\underline{\Delta}_{\wedge}$ satisfies the condition (E3) from \cite[Section 3.2]{SS}, i.e. its spectrum is contained in $(-\infty,0]$. 
\end{proof}

We close this section with a description of the domain of the complex powers of the model cone Laplacian.

\begin{corollary}
Let $c\geq0$, $z\in\mathbb{C}$ with $\mathrm{Re}(z)\in(0,1)$, $p\in(1,\infty)$, $\gamma$ be as in \eqref{gammaweight} and let $\underline{\Delta}_{\wedge}$ be the closed extension of $\Delta_{\wedge}$ given by \eqref{delta}. Then
$$
\mathcal{K}^{s+2\mathrm{Re}(z)+\varepsilon,\gamma+2\mathrm{Re}(z)+\varepsilon}_p(\partial\mathbb{B}^{\wedge})\oplus\mathbb{C}_{\omega}\hookrightarrow\mathcal{D}((c-\underline{\Delta}_{\wedge})^{z}) \hookrightarrow \mathcal{K}^{s+2\mathrm{Re}(z)-\varepsilon,\gamma+2\mathrm{Re}(z)-\varepsilon}_p(\partial\mathbb{B}^{\wedge})\oplus\mathbb{C}_{\omega},
$$
for every $\varepsilon>0$. 
\end{corollary}
\begin{proof}
The result follows by \cite[(I.2.5.2) and (I.2.9.6)]{Am}, \eqref{compare}, Lemma \ref{inter} and Corollary \ref{SecLap}.
\end{proof}

\section{The fractional Laplacian on a conic manifold}

The following result improves Theorem \ref{SecLapOrigin} from the sectoriality point of view and allows the definition of the fractional powers of the conic Laplacian.

\begin{theorem}\label{rsecdcon}
Let $p\in(1,\infty)$, $s\geq0$, $\gamma$ be as in \eqref{gammaweight} and let $\underline{\Delta}_{s}$ be the realization \eqref{delta1}. Then, $-\underline{\Delta}_{s}\in\mathcal{S}(\theta)$ for any $\theta\in[0,\pi)$.
\end{theorem}
\begin{proof}
By Theorem \ref{SecLapOrigin}, we know that $\sigma(\underline{\Delta}_{s})\subset (-\infty,0]$ and that $|\lambda|\|(\lambda-\underline{\Delta}_{s})^{-1}\|_{\mathcal{L}(X_{0}^{s})}$ is bounded in $\{\lambda\in S_{\theta}\, |\, |\lambda|\geq r\}$, for any $r>0$. Therefore we only have to check the behavior of $(\lambda-\underline{\Delta}_{s})^{-1}$ when $\lambda\rightarrow 0$, $\lambda\in S_{\theta}\backslash\{0\}$.

{\em Case of $s=0$ and $p=2$}. We follow the gluing construction ideas in \cite{Wm}, see also \cite[Section 3]{Ro2}. Let $\mathbb{M}$ be a closed connected smooth Riemannian manifold such that $(\mathcal{B}\backslash([0,1/2)\times\partial\mathcal{B}),g|_{\mathcal{B}\backslash([0,1/2)\times\partial\mathcal{B})})$ is isometrically embedded into $\mathbb{M}$. Denote by $L^{2}(\mathbb{M})$ the space of the square integrable functions on $\mathbb{M}$ with respect to the Riemannian measure. Let $\Delta_{\mathbb{M}}$ be the Laplacian on $\mathbb{M}$ and denote by $\underline{\Delta}_{\mathbb{M}}$ the unique self-adjoint extension of $\Delta_{\mathbb{M}}$ in $L^{2}(\mathbb{M})$. By the spectral theorem for self-adjoint operators, we have that $-\underline{\Delta}_{\mathbb{M}}\in\mathcal{S}(\theta)$.

Let $\omega_{1}$, $\omega_{2}$ be two cut-off functions with values on $[0,1]$ such that $\omega_{1}=\omega_{2}=1$ on $[0,1/2)\times\partial\mathcal{B}$, $\omega_{1}=\omega_{2}=0$ on $\mathcal{B}\backslash([0,1)\times\partial\mathcal{B})$ and $\omega_{2}=1$ on $\mathrm{supp}(\omega_{1})$. Moreover, let $\omega_{3}=1-\omega_{1}$ and let $\omega_{4}$ be a cut-off function such that $\omega_{4}=0$ on $[0,1/2)\times\partial\mathcal{B}$ and $\omega_{4}=1$ on $\mathrm{supp}(\omega_{3})$.

Consider the parametrix 
\begin{equation}\label{estQ}
Q(\lambda)=\omega_{1}(\lambda-\underline{\Delta}_{\wedge})^{-1}\omega_{2}+\omega_{3}(\lambda-\underline{\Delta}_{\mathbb{M}})^{-1}\omega_{4}, \quad \lambda\in\mathbb{C}\backslash(-\infty,0],
\end{equation}
where $\underline{\Delta}_{\wedge}$ is the closed extension \eqref{delta}. Clearly, $Q(\lambda)$ is a well defined map from $X_{0}^{0}$ to $X_{1}^{0}$, where the spaces $X_{0}^{0}$, $X_{1}^{0}$ are defined in Theorem \ref{SecLapOrigin}. Moreover
\begin{equation}\label{paramtrx}
(\lambda-\underline{\Delta}_{0})Q(\lambda)=I-P(\lambda), \quad \lambda\in\mathbb{C}\backslash(-\infty,0],
\end{equation}
where 
\begin{eqnarray}\nonumber
\lefteqn{P(\lambda)=[\underline{\Delta}_{\wedge},\omega_{1}](\lambda-\underline{\Delta}_{\wedge})^{-1}\omega_{2}+[\underline{\Delta}_{\mathbb{M}},\omega_{3}](\lambda-\underline{\Delta}_{\mathbb{M}})^{-1}\omega_{4}}\\\nonumber
&=&[\underline{\Delta}_{\wedge},\omega_{1}](c-\underline{\Delta}_{\wedge})^{-\eta}(c-\underline{\Delta}_{\wedge})^{\eta}(\lambda-\underline{\Delta}_{\wedge})^{-1}\omega_{2}\\\label{gtaet}
&&+[\underline{\Delta}_{\mathbb{M}},\omega_{3}](c-\underline{\Delta}_{\mathbb{M}})^{-\eta}(c-\underline{\Delta}_{\mathbb{M}})^{\eta}(\lambda-\underline{\Delta}_{\mathbb{M}})^{-1}\omega_{4},
\end{eqnarray}
for any $c>0$ and $\eta\in(1/2,1)$. Here the fractional powers of $c-\underline{\Delta}_{\mathbb{M}}$ can also be defined by the spectral theorem. 

In \eqref{gtaet} we regard the commutator $[\underline{\Delta}_{\wedge},\omega_{1}]$ as a first order cone differential operator on $\partial\mathbb{B}^{\wedge}$ and the commutator $[\underline{\Delta}_{\mathbb{M}},\omega_{3}]$ as first order differential operator on $\mathbb{M}$. Note that $[\underline{\Delta}_{\wedge},\omega_{1}](c-\underline{\Delta}_{\wedge})^{-\eta}$ and $[\underline{\Delta}_{\mathbb{M}},\omega_{3}](c-\underline{\Delta}_{\mathbb{M}})^{-\eta}$ are bounded operators that map to $\mathcal{H}_{2}^{\delta,\delta}(\mathbb{B})$ for certain $\delta>0$ sufficiently small. Since $\mathcal{H}_{2}^{\delta,\delta}(\mathbb{B})\hookrightarrow \mathcal{H}_{2}^{0,0}(\mathbb{B})$ is compact, see e.g. \cite[Remark 2.1 (b)]{Sei}, and similarly for the usual Sobolev spaces on $\mathbb{M}$, we have that $P(\lambda)$, $\lambda\in\mathbb{C}\backslash(-\infty,0]$, is a family of compact operators. Furthermore, since the residues of $(\lambda-\underline{\Delta}_{\wedge})^{-1}$ and $(\lambda-\underline{\Delta}_{\mathbb{M}})^{-1}$ are of finite rank (see e.g. \cite[Chapter III, Theorem 6.29]{Kato}), the residues of $P(\lambda)$ are of finite rank as well. Moreover, by the standard decay properties of the resolvent of a sectorial operator, see e.g. \cite[Lemma 2.3.3]{Tan}, we have that $\|P(\lambda)\|_{\mathcal{L}(X_{0}^{0})}\rightarrow0$ as $\lambda\rightarrow+\infty$. Therefore, by the meromorphic Fredholm theory, see e.g. \cite[Theorem XIII.13]{ReSi4}, $(I-P(\lambda))^{-1}$ exists outside a discrete set of points in $\mathbb{C}$. In combination with \cite[Theorem 4.1]{RS2}, we conclude that there exists a discrete set of points $D$ such that, due to \eqref{paramtrx}, we have
\begin{equation}\label{resolventconstr}
(\lambda-\underline{\Delta}_{0})^{-1}=Q(\lambda)(I-P(\lambda))^{-1}, \quad \lambda\in \mathbb{C}\backslash D.
\end{equation}

Next we argue that $0$ is a simple pole of $(\lambda-\underline{\Delta}_{0})^{-1}$. Let $\underline{\Delta}_{F}$ be the Friedrichs extension of $\Delta$ in $\mathcal{H}_{2}^{0,0}(\mathbb{B})$. By \cite[Section 5.3]{SS} its domain is given by 
$$
\mathcal{D}(\underline{\Delta}_{F})=\Bigg\{\begin{array}{lll} \mathcal{D}(\underline{\Delta}_{F,\min})\oplus\mathbb{C}_{\omega}\oplus\bigoplus_{q_{j}^{\pm}\in I_{0}, q_{j}^{\pm}<0}\mathcal{E}_{\Delta,0,q_{j}^{\pm}} & \text{if} & n=1\\
\mathcal{D}(\underline{\Delta}_{F,\min})\oplus\bigoplus_{q_{j}^{\pm}\in I_{0}, q_{j}^{\pm}\leq\frac{n-1}{2}}\mathcal{E}_{\Delta,0,q_{j}^{\pm}} & \text{if} & n\geq2,
\end{array}
$$ 
where $q_{j}^{\pm}$ and $\mathcal{E}_{\Delta,0,q_{j}^{\pm}}$ are described in \eqref{domflatcone}. Recall that $I_{0}=(\frac{n-3}{2},\frac{n+1}{2})$, so that when $n=2$, the pole $q_{j}^{\pm}=0$ is contained in $I_{0}$. Therefore, $\mathbb{C}_{\omega}\subset \mathcal{D}(\underline{\Delta}_{F})$ for $n=2$. If $n\geq3$, then $\mathbb{C}_{\omega}\subset \mathcal{D}(\underline{\Delta}_{F,\min})$ due to \eqref{dminlap}. We conclude that for any $n\geq1 $ we have
\begin{equation}\label{cembedtoDF}
\mathcal{D}(\underline{\Delta}_{0})\hookrightarrow \mathcal{D}(\underline{\Delta}_{F}).
\end{equation}

Assume first that $\gamma\geq0$. For any $\lambda\in\mathbb{C}\backslash(-\infty,0]$ according to \eqref{cembedtoDF} we have that
\begin{equation}\label{polecomp}
(\lambda-\underline{\Delta}_{F})^{-1}-(\lambda-\underline{\Delta}_{0})^{-1}=(\lambda-\underline{\Delta}_{F})^{-1}(\underline{\Delta}_{F}-\underline{\Delta}_{0})(\lambda-\underline{\Delta}_{0})^{-1}\quad \text{in}\quad \mathcal{H}_{2}^{0,\gamma}(\mathbb{B}).
\end{equation}
The right hand side of the above equation is $0$. Therefore,
$$
(\lambda-\underline{\Delta}_{F})^{-1}|_{\mathcal{H}_{2}^{0,\gamma}(\mathbb{B})}=(\lambda-\underline{\Delta}_{0})^{-1}
$$
when $\lambda\in\mathbb{C}\backslash(-\infty,0]$. This implies that the pole $0$ of $(\lambda-\underline{\Delta}_{0})^{-1}$ is simple.

Now assume that $\gamma<0$. The scalar product $\langle\cdot,\cdot\rangle$ of $\mathcal{H}_{2}^{0,0}(\mathbb{B})$ induces an identification of the dual space of $\mathcal{H}_{2}^{0,\gamma}(\mathbb{B})$ with $\mathcal{H}_{2}^{0,-\gamma}(\mathbb{B})$. The adjoint $\underline{\Delta}_{0}^{\ast}$ of $\underline{\Delta}_{0}$, is defined as usual by 
$$
\mathcal{D}(\underline{\Delta}_{0}^{\ast})=\big\{v\in \mathcal{H}_{2}^{0,-\gamma}(\mathbb{B})\, |\, \exists\, w\in \mathcal{H}_{2}^{0,-\gamma}(\mathbb{B}) \,\, \text{such that} \,\, \forall u\in \mathcal{D}(\underline{\Delta}_{0}) \,\, \text{we have} \,\, \langle v,\Delta u\rangle = \langle w,u\rangle\big\},
$$
see e.g. \cite[Section 5.3]{SS}. In particular, by \cite[Theorem 5.3]{SS} we have precisely
$$
\mathcal{D}(\underline{\Delta}_{0}^{\ast})=\bigg\{\begin{array}{lll} 
 \mathcal{H}_{2}^{2,-\gamma+2}(\mathbb{B})\oplus\mathbb{C}_{\omega} & \text{if} & n=1\\
 \mathcal{H}_{2}^{2,-\gamma+2}(\mathbb{B}) & \text{if} & n\geq2.
\end{array}
$$ 
We deduce that $\mathcal{D}(\underline{\Delta}_{0}^{\ast})\hookrightarrow \mathcal{D}(\underline{\Delta}_{F})$ and, similarly to \eqref{polecomp}, we have
$$
(\lambda-\underline{\Delta}_{F})^{-1}-(\lambda-\underline{\Delta}_{0}^{\ast})^{-1}=(\lambda-\underline{\Delta}_{F})^{-1}(\underline{\Delta}_{F}-\underline{\Delta}_{0}^{\ast})(\lambda-\underline{\Delta}_{0}^{\ast})^{-1}\quad \text{in}\quad \mathcal{H}_{2}^{0,-\gamma}(\mathbb{B}).
$$
Again, the right hand side of the above equation is identical to $0$ and we conclude that $0$ is a simple pole of $(\lambda-\underline{\Delta}_{0}^{\ast})^{-1}$. By \cite[Proposition 1.3 (v)]{DHP} we find that $0$ is a simple pole of $(\lambda-\underline{\Delta}_{0})^{-1}$ as well. 

{\em Case of $s=0$ and $p\in(1,\infty)$}. Denote by $R_{p}(\lambda)$ the resolvent $(\lambda-\underline{\Delta}_{0})^{-1}$ in the space $\mathcal{H}_{p}^{0,\gamma}(\mathbb{B})$. In \cite[Proposition 3.1]{RS1} it has been shown that $\underline{\Delta}_{0}$ satisfies the ellipticity conditions (E1), (E2) and (E3) of \cite[Section 3.2]{SS}. Therefore, by \cite[Theorem 4.1]{SS} for each $p\in(1,\infty)$ there exists some $r_{0}>0$ such that $R_{p}(\lambda)$ exists for $\lambda\in S_{\theta}$, $|\lambda|\geq r_{0}$, and is equal to $R_{2}(\lambda)$, in the sense that $R_{p}(\lambda)$ is the restriction of $R_{2}(\lambda)$ and vice versa. Furthermore, by \cite[Theorem 4.2]{RS2} we know that $R_{p}(\lambda)$ exists for all $\lambda\in S_{\theta}\backslash\{0\}$ and that, for each fixed $\delta_{0}>0$, $\|R_{p}(\lambda)\|_{\mathcal{L}(\mathcal{H}_{p}^{0,\gamma}(\mathbb{B}))}$ is uniformly bounded by $K_{0}$ when $\lambda\in S_{\theta}$, $|\lambda|\geq\delta_{0}$, for certain $K_{0}>0$. Thus, 
\begin{equation}\label{opsdt}
\|R_{2}(\lambda_{0})\|_{\mathcal{L}(\mathcal{H}_{p}^{0,\gamma}(\mathbb{B}))}\leq K_{0} \quad \text{when} \quad \lambda_{0}\in S_{\theta}, \, |\lambda_{0}|\geq r_{0}, 
\end{equation}
and by Neumann series we get that
$$
R_{2}(\lambda)=R_{2}(\lambda_{0})\sum_{k=0}^{\infty}((\lambda_{0}-\lambda)R_{2}(\lambda_{0}))^{k} \in \mathcal{L}(\mathcal{H}_{p}^{0,\gamma}(\mathbb{B})), \quad \lambda,\lambda_{0}\in S_{\theta}, \, |\lambda-\lambda_{0}|\leq \frac{1}{2K_{0}}, \, |\lambda_{0}|=r_{0}.
$$
Therefore, by analyticity, i.e. by the identity theorem, we obtain that $R_{p}(\lambda)=R_{2}(\lambda)$ for all $\lambda\in S_{\theta}$ with $|\lambda|\geq r_{0}-(2K_{0})^{-1}$. Then, \eqref{opsdt} holds for $\lambda_{0}\in S_{\theta}$, $|\lambda_{0}|\geq r_{0}-(2K_{0})^{-1}$, with the same bound $K_{0}$. After finitely many steps we deduce that $R_{p}(\lambda)=R_{2}(\lambda)$ for all $\lambda\in S_{\theta}$, $|\lambda|\geq 2\delta_{0}$, and the result follows since $\delta_{0}>0$ was arbitrary.

{\em Case of $s>0$ and $p\in(1,\infty)$}. In the Step 1 of the proof of \cite[Theorem 3.3]{RS2} we have seen that the resolvent of $\underline{\Delta}_{s}$ is the restriction of the resolvent of $\underline{\Delta}_{0}$ to $\mathcal{H}_{p}^{s,\gamma}(\mathbb{B})$. Moreover, let $B_{-1}$ be the residue of $(\lambda-\underline{\Delta}_{0})^{-1}$ at $\lambda=0$. 

Assume first that $s\in[0,2]$. By Lemma \ref{expandatzero} we have that 
$$
B_{-1}\in\mathcal{L}(\mathcal{H}_{p}^{0,\gamma}(\mathbb{B}),\mathcal{D}(\underline{\Delta}_{0}))\hookrightarrow \mathcal{L}(\mathcal{H}_{p}^{0,\gamma}(\mathbb{B}),\mathcal{H}_{p}^{2,\gamma}(\mathbb{B})).
$$
Therefore, $B_{-1}\in\mathcal{L}(\mathcal{H}_{p}^{s,\gamma}(\mathbb{B}))$. Hence, the result follows by the previous result for the case of $s=0$ and by \cite[Theorem 4.1]{RS2}. 

Next, let that $s\in[2,4]$. Lemma \ref{expandatzero} implies that 
$$
B_{-1}\in\mathcal{L}(\mathcal{H}_{p}^{2,\gamma}(\mathbb{B}),\mathcal{D}(\underline{\Delta}_{2}))\hookrightarrow \mathcal{L}(\mathcal{H}_{p}^{2,\gamma}(\mathbb{B}),\mathcal{H}_{p}^{4,\gamma}(\mathbb{B})).
$$
Hence, we have that $B_{-1}\in\mathcal{L}(\mathcal{H}_{p}^{s,\gamma}(\mathbb{B}))$ and the result follows by the result for the case of $s=2$ and by \cite[Theorem 4.1]{RS2}. Iteration then shows the assertion.
\end{proof}

We are now in a position to prove the main result of this section.

\subsubsection*{Proof of Theorem 1.2}
The fractional Laplacian $(-\underline{\Delta}_{s})^{\sigma}$ is defined by Theorem \ref{fracpower} and Theorem \ref{rsecdcon}. We denote 
\begin{equation}\label{fracdomD}
X_{\sigma}^{s}=\mathcal{D}((-\underline{\Delta}_{s})^{\sigma}). 
\end{equation}
By \eqref{compare}, $X_{\sigma}^{s}=\mathcal{D}((\delta-\underline{\Delta}_{s})^{\sigma})$ for any $\delta>0$, so that the embedding \eqref{domfrac} follows by \cite[Corollary 5.3]{RS3}. Concerning the sharp description of the domain \eqref{sharpdom} under \eqref{gammapoles}, it follows by \cite[Lemma 4.5]{LR} and \eqref{compare}. The $R$-sectoriality for the fractional Laplacian follows by Theorem \ref{Rsec}, Theorem \ref{SecLapOrigin} and Theorem \ref{rsecdcon}. \mbox{\ } \hfill $\square$

\section{The fractional porous medium equation}

The starting point for the study of our fractional diffusion is the following observation, which shows that the commutator between a function of certain regularity and the fractional Laplacian is of lower order in a fractional sense.

\begin{lemma}[Commutation]\label{commlemma}
Let $p\in(1,\infty)$, $s\geq0$, $\gamma$ be as in \eqref{gammaweight}, $c>0$, $\sigma\in(0,1)$, $\eta\in(\frac{1}{2},1)$, $\mu>s+1+2\eta+\frac{n+1}{p}$, $\xi>\max\{\gamma+2,\frac{n+3}{2}\}$, $w\in \mathcal{H}_{p}^{\mu,\xi}(\mathbb{B})\oplus\mathbb{C}_{\omega}$ and
$$
\rho_{0}=\bigg\{\begin{array}{lll} \eta-\frac{1}{2} &\text{if} & \xi\geq\gamma+2\eta+1\\ \frac{\xi-\gamma}{2}-1 & \text{if} & \xi<\gamma+2\eta+1.\end{array}
$$ 
Then, for each $\nu>\sigma+\eta-1$ and $\rho\in[0,\rho_{0})$ we have
$$
[w,(c-\underline{\Delta}_{s})^{\sigma}]\in \mathcal{L}(\mathcal{D}((c-\underline{\Delta}_{s})^{\nu}),\mathcal{D}((c-\underline{\Delta}_{s})^{\rho})).
$$
Consequently, the above commutator is of lower order as we can choose $\nu<\sigma$, and, on the other hand, it maps to the domain of the $\rho$-power. 
\end{lemma}
\begin{proof}
In local coordinates $(x,y)\in(0,1)\times\partial\mathcal{B}$ on the collar part the first order differential operator $[\underline{\Delta}_{s},w]$ is of the form 
$$
[\underline{\Delta}_{s},w]=2(\partial_{x}w)\partial_{x}+2x^{-2}\langle \nabla w,\nabla \cdot\rangle_{h}+\big((\partial_{x}^{2}w)+nx^{-1}(\partial_{x}w)+x^{-2}\Delta_{h}w\big),
$$
where $\langle \cdot,\cdot\rangle_{h}$ and $\nabla$ are respectively the Riemannian scalar product and the gradient on $\partial\mathbb{B}$. Thus, if we denote $A_{s}=c-\underline{\Delta}_{s}$ and, according to \eqref{domfrac}, write any $u\in \mathcal{D}(A_{s}^{\eta})$ as $u=u_{\mathcal{H}}+u_{\mathbb{C}}$ with 
$$
u_{\mathcal{H}}\in \bigcap_{\varepsilon>0} \mathcal{H}_{p}^{s+2\eta-\varepsilon,\gamma+2\eta-\varepsilon}(\mathbb{B})\quad \text{and}\quad u_{\mathbb{C}}\in\mathbb{C}_{\omega},
$$
then by the regularity of $w$ and \cite[Corollary 3.3]{RS3}, we conclude that 
$$
[\underline{\Delta}_{s},w]A_{s}^{-\eta}\in\bigcap_{\varepsilon>0}\mathcal{L}(X_{0}^{s},\mathcal{H}_{p}^{s+2\eta-1-\varepsilon,\tau_{\varepsilon}}(\mathbb{B})), \quad \text{where} \quad \tau_{\varepsilon}=\min\{\gamma+2\eta-1-\varepsilon,\xi-2\}.
$$ 

Therefore, by \eqref{Aexp}, in $\mathcal{D}(\underline{\Delta}_{s})$ we have
\begin{eqnarray*}
w A_{s}^{\sigma}-A_{s}^{\sigma}w&=&\frac{\sin(\pi \sigma)}{\pi}\int_{0}^{\infty}x^{\sigma-1}([\underline{\Delta}_{s},w](A_{s}+x)^{-1}+A_{s}[w,(A_{s}+x)^{-1}])dx\\
&=&\frac{\sin(\pi \sigma)}{\pi}\int_{0}^{\infty}x^{\sigma-1}\Big([\underline{\Delta}_{s},w]A_{s}^{-\eta}A_{s}^{\eta-\nu}(A_{s}+x)^{-1}\\
&&+A_{s}(A_{s}+x)^{-1}[w,\underline{\Delta}_{s}]A_{s}^{-\eta}A_{s}^{\eta-\nu}(A_{s}+x)^{-1}\Big)A_{s}^{\nu}dx.
\end{eqnarray*}
Due to \cite[Lemma 2.3.3]{Tan} (alternatively see Lemma \ref{lemmadecay} in the Appendix), the right hand side of the above equation belongs to $\mathcal{L}(\mathcal{D}(A_{s}^{\nu}),X_{0}^{s})$, so that $[w,A_{s}^{\sigma}]\in \mathcal{L}(\mathcal{D}(A_{s}^{\nu}),X_{0}^{s})$.

Moreover, by \eqref{domfrac}, the integral
\begin{eqnarray*}
\lefteqn{\int_{0}^{\infty}x^{\sigma-1}\Big(A_{s}^{\rho}[\underline{\Delta}_{s},w]A_{s}^{-\eta}A_{s}^{\eta-\nu}(A_{s}+x)^{-1}}\\
&&+A_{s}(A_{s}+x)^{-1}A_{s}^{\rho}[w,\underline{\Delta}_{s}]A_{s}^{-\eta}A_{s}^{\eta-\nu}(A_{s}+x)^{-1}\Big)A_{s}^{\nu}dx
\end{eqnarray*}
converges absolutely, which implies that $[w,A_{s}^{\sigma}]\in \mathcal{L}(\mathcal{D}(A_{s}^{\nu}),\mathcal{D}(A_{s}^{\rho}))$. 
\end{proof}

Next we show $R$-sectoriality (and hence maximal $L^{q}$-regularity) for the linearization of \eqref{FPME1}. The above commutation property allows us to extend the freezing-of-coefficients method to our non-local situation. The resulting method is applicable to the more general case of linear combinations of terms, each one being a product of a function and a fractional power of a local operator.

\begin{theorem}\label{wperturb}
Let $p\in(1,\infty)$, $s=0$, $\gamma$ be chosen as in \eqref{gammaweight}, $\sigma\in(0,1)$ and $(-\underline{\Delta}_{0})^{\sigma}$ be the fractional Laplacian defined in Theorem \ref{ThFracLap}. If 
$$
w\in \bigcup_{\varepsilon>0}\mathcal{H}_{p}^{\frac{n+1}{p}+\varepsilon,\frac{n+1}{2}+\varepsilon}(\mathbb{B})\oplus\mathbb{C}_{\omega}
$$ 
satisfies $w\geq \alpha>0$ on $\mathbb{B}$, for certain $\alpha>0$, then for each $\theta\in[0,\pi)$ there exists a $c>0$ such that $w(-\underline{\Delta}_{0})^{\sigma}+c\in \mathcal{R}(\theta)$.
\end{theorem}
\begin{proof}
{\em Step 1: $R$-sectoriality for $w(c_{0}-\underline{\Delta}_{0})^{\sigma}+c$}.
By \cite[Lemma 3.2]{RS3} we have that
$$
\bigcup_{\varepsilon>0}\mathcal{H}_{p}^{\frac{n+1}{p}+\varepsilon,\frac{n+1}{2}+\varepsilon}(\mathbb{B})\oplus\mathbb{C}_{\omega} \hookrightarrow C(\mathbb{B}).
$$
Let $c_{0}>0$ and denote $A=c_{0}-\underline{\Delta}_{0}$. By the identity 
$$
\lambda(w(z_{0})A^{\sigma}+\lambda)^{-1}=(\lambda/w(z_{0}))(A^{\sigma}+\lambda/w(z_{0}))^{-1}, \quad \lambda\in S_{\theta}\backslash\{0\}, \, z_{0}\in\mathbb{B}, 
$$
we deduce that $w(z_{0})A^{\sigma}$ is $R$-sectorial and its $R$-sectorial bound is uniformly bounded in $z_{0}\in\mathbb{B}$ by the $R$-sectorial bound of $A^{\sigma}$. Hence, due to \cite[Lemma 2.6]{RS3}, the $R$-sectorial bound of $w(z_{0})A^{\sigma}+c$ is uniformly bounded in $z_{0}\in\mathbb{B}$ and $c>0$.

Let $r>0$ and choose an open cover of $\mathbb{B}$ consisting of balls $B_{j}=B_{r}(z_{j})$, $z_{j}\in\mathbb{B}^\circ$, $j\in\{1,\dots,N\}$, of radius $r$, together with a collar neighborhood $B_{0}= [0,r)\times\partial \mathcal{B}$. We assume that $\overline{B_{3r/2}(z_{j})}$, $j\in\{1,\dots,N\}$, do not intersect $\{0\}\times\partial\mathcal{B}$. Let $\widetilde{\omega} :\mathbb{R}\rightarrow [0,1]$ be a smooth non-increasing function that equals $1$ on $[0,1/2]$ and $0$ on $[3/4,\infty)$, and denote by $d=d(z,\widetilde{z})$ the geodesic distance between two points $z,\widetilde{z}\in\mathbb{B}$ with respect to the metric $g$. Fix some $z_{0}\in \{0\}\times\partial\mathcal{B}$ and define 
$$
w_{j}(z)=\widetilde{\omega}\Big(\frac{d(z,z_j)}{2r}\Big)w(z) + \Big(1-\widetilde{\omega}\Big(\frac{d(z,z_j)}{2r}\Big)\Big)w(z_j), \quad z\in\mathbb{B}, \quad j\in\{0,\dots,N\}.
$$
Since $\|w(z_{j})-w_{j}(\cdot)\|_{C(\mathbb{B})}$, and therefore the norm of $w(z_{j})-w_{j}(\cdot)$ as a multiplier on $X_{0}^{0}$, becomes arbitrarily small when $r\rightarrow 0$ (for $j=0$ recall that by \cite[Lemma 2.6]{RS3} $w$ is constant along the boundary), by writing
$$
w_{j}A^{\sigma}+c=w(z_{j})A^{\sigma}+c+(w_{j}-w(z_{j}))A^{\sigma},
$$
from \cite[Proposition 4.4.2]{PS} we see that for small values of $r$ each $w_{j}A^{\sigma}+c$ becomes $R$-sectorial of angle $\theta$. 

Moreover, by interpolation, see e.g. \cite[(I.2.5.2), (I.2.8.4) and (I.2.9.6)]{Am}, for each $0<\xi_{0}<\xi_{1}<\xi_{2}<1$ and $j\in\{0,\dots,N\}$ we have
\begin{equation}\label{cpomparecomplex}
\mathcal{D}(A^{\sigma\xi_{2}})\hookrightarrow \mathcal{D}((w_{j}A^{\sigma}+c)^{\xi_{1}})\hookrightarrow \mathcal{D}(A^{\sigma\xi_{0}}).
\end{equation}

{\em Left inverse.} Let $\phi_{j}\in C^{\infty}(\mathbb{B})$, $j=\{0,\dots,N\}$, be a partition of unity that is subordinate to $\{B_{j}\}_{ j\in\{0,\dots,N\}}$ and let $\psi_{j}\in C^{\infty}(\mathbb{B})$, $ j\in\{0,\dots,N\}$, with values on $[0,1]$ and supported in $B_{j}$ such that $\psi_{j}=1$ on $\mathrm{supp}(\phi_{j})$. Recall the notation \eqref{fracdomD}. If $u\in X_{\sigma}^{0}$, $f\in X_{0}^{0}$ and $\lambda\in S_{\theta}$, then by multiplying
$$
(wA^{\sigma}+c+\lambda)u=f
$$
with $\phi_{j}$, $j=\{0,\dots,N\}$, and notting that $\phi_{j}w=\phi_{j}w_{j}$, we obtain 
$$
\phi_{j}w_{j}A^{\sigma}u+(c+\lambda)\phi_{j}u=\phi_{j}f,
$$
and hence
$$
(w_{j}A^{\sigma}+c+\lambda) \phi_{j}u=\phi_{j}f+[w_{j}A^{\sigma},\phi_{j}]u.
$$
By applying the resolvent of $w_{j}A^{\sigma}+c+\lambda$ to the above equation we get
$$
\phi_{j}u=(w_{j}A^{\sigma}+c+\lambda)^{-1}(\phi_{j}f+w_{j}[A^{\sigma},\phi_{j}]u),
$$
where by multiplying with $\psi_{j}$ and then summing up we obtain 
\begin{equation}\label{leftinv}
u=\sum_{j=0}^{N}\psi_{j}(w_{j}A^{\sigma}+c+\lambda)^{-1}\phi_{j}f+\sum_{j=0}^{N}\psi_{j}(w_{j}A^{\sigma}+c+\lambda)^{-1}w_{j}[A^{\sigma},\phi_{j}]u.
\end{equation}

Fix $\eta>0$ such that $w\in\mathcal{H}_{p}^{\frac{n+1}{p}+\eta,\frac{n+1}{2}+\eta}(\mathbb{B})\oplus\mathbb{C}_{\omega}$. Due to Lemma \ref{commlemma}, we write 
$$
(w_{j}A^{\sigma}+c+\lambda)^{-1}w_{j}[A^{\sigma},\phi_{j}]=(w_{j}A^{\sigma}+c+\lambda)^{-1}w_{j}(c_{0}-\underline{\Delta}_{0})^{-\rho}(c_{0}-\underline{\Delta}_{0})^{\rho}[A^{\sigma},\phi_{j}],
$$
for some $\rho\in (0,\frac{\eta}{2})$, so that $(c_{0}-\underline{\Delta}_{0})^{\rho}[A^{\sigma},\phi_{j}]\in \mathcal{L}(X_{\sigma}^{0},X_{0}^{0})$. In addition, by \eqref{domfrac} and \eqref{compare} 
$$
\mathcal{D}((c_{0}-\underline{\Delta}_{0})^{\rho})\hookrightarrow \bigcap_{\varepsilon>0}\mathcal{H}_{p}^{2\rho-\varepsilon,\gamma+2\rho-\varepsilon}(\mathbb{B})\oplus\mathbb{C}_{\omega}.
$$
Moreover, by \cite[Lemma 3.3]{RS3} each $w_{j}$ acts by multiplication as a bounded map on 
$$
\bigcap_{\varepsilon>0}\mathcal{H}_{p}^{2\rho-\varepsilon,\gamma+2\rho-\varepsilon}(\mathbb{B})\oplus\mathbb{C}_{\omega}.
$$
Therefore, by \eqref{domfrac}, \eqref{compare} and \eqref{cpomparecomplex} we obtain that
$$
w_{j}(c_{0}-\underline{\Delta}_{0})^{-\rho}\in \mathcal{L}(X_{0}^{0},\mathcal{D}((w_{j}A^{\sigma}+c)^{\widetilde{\rho}})),
$$ 
for certain $\widetilde{\rho}\in (0,\rho)$. Hence, from \cite[Lemma 2.3.3]{Tan}, by taking $c>0$ sufficiently large, the $\mathcal{L}(X_{\sigma}^{0})$ norm of the second term on the right hand side of \eqref{leftinv} becomes arbitrary small, uniformly in $\lambda\in S_{\theta}$. We conclude that there exists some $\widetilde{c}>0$ such that for $c\geq\widetilde{c}$ the operator $wA^{\sigma}+c+\lambda$ has a left inverse $L$ that belongs to $\mathcal{L}(X_{0}^{0},X_{\sigma}^{0})$; in particular
\begin{equation}\label{leftinvexp}
L=\sum_{k=0}^{\infty}Q^{k}(\lambda)R(\lambda),
\end{equation}
where 
$$
Q(\lambda)=\sum_{j=0}^{N}\psi_{j}(w_{j}A^{\sigma}+c+\lambda)^{-1}w_{j}[A^{\sigma},\phi_{j}] \quad \text{and} \quad R(\lambda)=\sum_{j=0}^{N}\psi_{j}(w_{j}A^{\sigma}+c+\lambda)^{-1}\phi_{j}.
$$ 

{\em Right inverse.} By notting that $\psi_{j}w=\psi_{j}w_{j}$, $j=\{0,\dots,N\}$, from \eqref{leftinv} we obtain 
\begin{eqnarray}\nonumber
\lefteqn{(wA^{\sigma}+c+\lambda)L=I+w\sum_{j=0}^{N}[A^{\sigma},\psi_{j}](w_{j}A^{\sigma}+c+\lambda)^{-1}\phi_{j}}\\\label{rightinv}
&&+\sum_{j=0}^{N}\psi_{j}w_{j}[A^{\sigma},\phi_{j}]L+w\sum_{j=0}^{N}[A^{\sigma},\psi_{j}](w_{j}A^{\sigma}+c+\lambda)^{-1}w_{j}[A^{\sigma},\phi_{j}]L,
\end{eqnarray}
where we have used the fact that $\psi_{j}\phi_{j}=\phi_{j}$ and $\sum_{j=1}^{N}\phi_{j}=1$. 

Let $\nu\in(\sigma,1)$ and write 
$$
[A^{\sigma},\phi_{i}]\psi_{j}(w_{j}A^{\sigma}+c+\lambda)^{-1}=[A^{\sigma},\phi_{i}]\psi_{j}(w_{j}A^{\sigma}+\widetilde{c})^{-\nu}(w_{j}A^{\sigma}+\widetilde{c})^{\nu}(w_{j}A^{\sigma}+c+\lambda)^{-1};
$$ 
here $[A^{\sigma},\phi_{i}]\psi_{j}(w_{j}A^{\sigma}+\widetilde{c})^{-\nu}\in\mathcal{L}(X_{0}^{0})$ for all $i,j\in\{0,\dots,N\}$ due to Lemma \ref{commlemma} and \eqref{cpomparecomplex}. Therefore, by writing 
$$
Q^{k}(\lambda)=Q(\lambda)(c_{0}-\underline{\Delta}_{0})^{-\nu}(c_{0}-\underline{\Delta}_{0})^{\nu}\cdots(c_{0}-\underline{\Delta}_{0})^{-\nu}(c_{0}-\underline{\Delta}_{0})^{\nu}Q(\lambda), \quad k\in\mathbb{N}_{0},
$$
in \eqref{leftinvexp} and using \cite[Lemma 2.3.3]{Tan}, we see that $\|[A^{\sigma},\phi_{j}]L\|_{\mathcal{L}(X_{0}^{0})}$, $j\in\{0,\dots,N\}$, becomes arbitrary small, uniformly in $\lambda\in S_{\theta}$, by taking $c\geq\widetilde{c}$ sufficiently large. 

Similarly, we write
$$
[A^{\sigma},\psi_{j}](w_{j}A^{\sigma}+c+\lambda)^{-1}=[A^{\sigma},\psi_{j}](w_{j}A^{\sigma}+\widetilde{c})^{-\nu}(w_{j}A^{\sigma}+\widetilde{c})^{\nu}(w_{j}A^{\sigma}+c+\lambda)^{-1}, \quad j\in\{0,\dots,N\},
$$ 
so that $[A^{\sigma},\psi_{j}](w_{j}A^{\sigma}+\widetilde{c})^{-\nu}\in\mathcal{L}(X_{0}^{0})$ due to Lemma \ref{commlemma} and \eqref{cpomparecomplex}. Hence, by \cite[Lemma 2.3.3]{Tan} the last three terms on the right hand side of \eqref{rightinv} become arbitrary small, uniformly in $\lambda\in S_{\theta}$, by taking $c\geq\widetilde{c}$ sufficiently large. This provides us for large $c>0$ a right inverse for $wA^{\sigma}+c+\lambda$ which belongs to $\mathcal{L}(X_{0}^{0},X_{\sigma}^{0})$.

{\em $R$-sectoriality.} Denote by $K\geq1$ the maximum of all $R$-sectorial bounds of $w_{j}A^{\sigma}+c$, $j\in\{0,\dots,N\}$; recall that, due to \cite[Lemma 2.6]{RS3}, $K$ can be chosen independent of $c\geq \widetilde{c}$. Let $\lambda_{1},\dots,\lambda_{M}\in S_{\theta}\backslash\{0\}$, $M\in\mathbb{N}$, $v_{1},\dots,v_{M}\in X_{0}^{0}$ and $\{\epsilon_{k}\}_{k\in\mathbb{N}}$ be the sequence of the Rademacher functions. We have that
\begin{eqnarray}\nonumber
\|\sum_{i=1}^{M}\varepsilon_{i}\lambda_{i}R(\lambda_{i})v_{i}\|_{L^{2}(0,1;X_{0}^{0})}&\leq&\sum_{j=0}^{N}\sup(|\psi_{j}|)\|\sum_{i=1}^{M}\varepsilon_{i}\lambda_{i}(w_{j}A^{\sigma}+c+\lambda)^{-1}\phi_{j}v_{i}\|_{L^{2}(0,1;X_{0}^{0})}\\\label{Rest}
&\leq&K(N+1)\|\sum_{i=1}^{M}\varepsilon_{i}v_{i}\|_{L^{2}(0,1;X_{0}^{0})}.
\end{eqnarray} 

Moreover, due to 
$$
[A^{\sigma},\phi_{j}]\psi_{k}(w_{k}A^{\sigma}+c)^{-1}=[A^{\sigma},\phi_{j}]\psi_{k}(w_{k}A^{\sigma}+\widetilde{c})^{-\nu}(w_{k}A^{\sigma}+\widetilde{c})^{\nu}(w_{k}A^{\sigma}+c)^{-1}, \quad j,k\in\{0,\dots,N\},
$$
\cite[Lemma 2.3.3]{Tan} and \eqref{cpomparecomplex}, for each $\delta>0$ there exists a $c\geq\widetilde{c}$ such that 
$$
\|[A^{\sigma},\phi_{j}]\psi_{k}(w_{k}A^{\sigma}+c)^{-1}\|_{\mathcal{L}(X_{0}^{0})}<\delta. 
$$
Hence, for each $l\in \mathbb{N}$ we estimate 
\begin{eqnarray*}
\lefteqn{\|\sum_{i=1}^{M}\varepsilon_{i}\lambda_{i}Q^{l}(\lambda_{i})R(\lambda_{i})v_{i}\|_{L^{2}(0,1;X_{0}^{0})}}\\
&\leq&\sum_{j=0}^{N}\sup(|\psi_{j}|)\|\sum_{i=1}^{M}\varepsilon_{i}\lambda_{i}(w_{j}A^{\sigma}+c+\lambda_{i})^{-1}w_{j}[A^{\sigma},\phi_{j}]Q^{l-1}(\lambda_{i})R(\lambda_{i})v_{i}\|_{L^{2}(0,1;X_{0}^{0})}\\
&\leq&K(N+1)\max_{j}\|\sum_{i=1}^{M}\varepsilon_{i}[A^{\sigma},\phi_{j}]\\
&&\Big(\sum_{k=0}^{N}\psi_{k}(w_{k}A^{\sigma}+c+\lambda_{i})^{-1}w_{k}[A^{\sigma},\phi_{k}]\Big)Q^{l-2}(\lambda_{i})R(\lambda_{i})v_{i}\|_{L^{2}(0,1;X_{0}^{0})}\\
&\leq&\delta K(K+1)(N+1)^{2}\max_{k}\|\sum_{i=1}^{M}\varepsilon_{i}[A^{\sigma},\phi_{k}]Q^{l-2}(\lambda_{i})R(\lambda_{i})v_{i}\|_{L^{2}(0,1;X_{0}^{0})}\\
&\leq&\delta^{l-1} ((N+1)(K+1))^{l}\max_{k}\|\sum_{i=1}^{M}\varepsilon_{i}[A^{\sigma},\phi_{k}]R(\lambda_{i})v_{i}\|_{L^{2}(0,1;X_{0}^{0})}\\
&\leq&\delta^{l}(N+1)((N+1)(K+1))^{l+1}\|\sum_{i=1}^{M}\varepsilon_{i}v_{i}\|_{L^{2}(0,1;X_{0}^{0})}.
\end{eqnarray*} 
By taking $\delta<(2(N+1)(K+1))^{-1}$, from \eqref{leftinvexp}, \eqref{Rest} and the above inequality we conclude that, for $c>0$ sufficiently large, $wA^{\sigma}+c$ is $R$-sectorial and its $R$-sectorial bound is bounded by $(N+1)K+(N+1)^{2}(K+1)$.

{\em Step 2: $R$-sectoriality for $w(-\underline{\Delta}_{0})^{\sigma}+c$}. Let $c_{1}>0$ be fixed and sufficiently large. By \cite[Lemma 2.6]{RS3} and the estimate in the part (i) in the proof of Theorem \ref{Rsec}, the $R$-sectorial bound of $(c-\underline{\Delta}_{0})^{\sigma}$ is uniformly bounded in $c\geq c_{1}$. By the Step 1 above and \cite[Lemma 2.6]{RS3}, both operators $w(c-\underline{\Delta}_{0})^{\sigma}+c_{1}$ and $w(c-\underline{\Delta}_{0})^{\sigma}+c^{\sigma+\xi}$ are $R$-sectorial and their $R$-sectorial bounds are uniformly bounded in $c\geq c_{1}$, where $\xi>1$ is fixed. By \eqref{compare} we estimate
\begin{eqnarray*}
\lefteqn{\|(w(c-\underline{\Delta}_{0})^{\sigma}-w(-\underline{\Delta}_{0})^{\sigma})(w(c-\underline{\Delta}_{0})^{\sigma}+c^{\sigma+\xi})^{-1}\|_{\mathcal{L}(X_{0}^{0})}}\\
&\leq&\|w\|_{\mathcal{L}(X_{0}^{0})}\|(c-\underline{\Delta}_{0})^{\sigma}-(-\underline{\Delta}_{0})^{\sigma})((c-\underline{\Delta}_{0})^{\sigma}+c^{\sigma+\xi})^{-1}\|_{\mathcal{L}(X_{0}^{0})}\\
&&\times\|((c-\underline{\Delta}_{0})^{\sigma}+c^{\sigma+\xi})(w(c-\underline{\Delta}_{0})^{\sigma}+c^{\sigma+\xi})^{-1}\|_{\mathcal{L}(X_{0}^{0})}\\
&\leq&C_{0}\|w\|_{\mathcal{H}_{p}^{\frac{n+1}{p}+\eta,\frac{n+1}{2}+\eta}(\mathbb{B})\oplus\mathbb{C}_{\omega}}M_{0} c^{\sigma}\frac{K_{0}}{c^{\sigma+\xi}}\\
&&\times\|w^{-1}\|_{\mathcal{L}(X_{0}^{0})}\|w(c-\underline{\Delta}_{0})^{\sigma}+c^{\sigma+\xi}+(w-1)c^{\sigma+\xi})(w(c-\underline{\Delta}_{0})^{\sigma}+c_{1}+c^{\sigma+\xi}-c_{1})^{-1}\|_{\mathcal{L}(X_{0}^{0})}\\
&\leq&c^{-\xi}C_{1}\|w\|_{\mathcal{H}_{p}^{\frac{n+1}{p}+\eta,\frac{n+1}{2}+\eta}(\mathbb{B})\oplus\mathbb{C}_{\omega}}\|w^{-1}\|_{\mathcal{H}_{p}^{\frac{n+1}{p}+\eta,\frac{n+1}{2}+\eta}(\mathbb{B})\oplus\mathbb{C}_{\omega}}\\
&&\times \Big(1+\|w-1\|_{\mathcal{H}_{p}^{\frac{n+1}{p}+\eta,\frac{n+1}{2}+\eta}(\mathbb{B})\oplus\mathbb{C}_{\omega}}c^{\sigma+\xi}\frac{\widetilde{K}_{0}}{1+c^{\sigma+\xi}-c_{1}}\Big),
\end{eqnarray*}
for certain $C_{0},C_{1}, M_{0} >0$, where $K_{0}$ is the sectorial bound of $(c-\underline{\Delta}_{0})^{\sigma}\in \mathcal{S}(0)$ and $\widetilde{K}_{0}$ is the sectorial bound of $w(c-\underline{\Delta}_{0})^{\sigma}+c_{1}\in \mathcal{P}(0)$. By taking $c\geq c_{1}$ sufficiently large, we obtain the result by perturbation (see \cite[Proposition 4.4.2]{PS}) due to
$$
w(-\underline{\Delta}_{0})^{\sigma}+c^{\sigma+\xi}=w(c-\underline{\Delta}_{0})^{\sigma}+c^{\sigma+\xi}+w(-\underline{\Delta}_{0})^{\sigma}-w(c-\underline{\Delta}_{0})^{\sigma}.
$$
\end{proof}

If the multiplication function has better regularity, then, due to Theorem \ref{rsecdcon}, we expect to have a result similar to Theorem \ref{wperturb} with $(-\underline{\Delta}_{0})^{\sigma}$ replaced by $(-\underline{\Delta}_{s})^{\sigma}$, $s>0$. However, by following the proof of Theorem \ref{wperturb}, due to \cite[Corollary 3.3]{RS3}, the norm of $w(z_{j})-w_{j}(\cdot)$ as a multiplier on $X_{0}^{s}$ is determined by $\|w(z_{j})-w_{j}(\cdot)\|_{Y}$ and not by $\|w(z_{j})-w_{j}(\cdot)\|_{C(\mathbb{B})}$ anymore, where $Y$ is a Mellin-Sobolev space of certain order that depends on $s$. Therefore, it cannot become arbitrarily small when $r\rightarrow 0$. On the other hand, by Theorem \eqref{HaHiTh}, the $R$-sectoriality result of Theorem \ref{wperturb} can be extended to higher order Mellin-Sobolev spaces as follows.

\begin{theorem}\label{wperturbsect}
Let $\sigma\in(0,1)$, $p\in(1,\infty)$, $s\geq0$, $\gamma$ be as \eqref{gammaweight} and \eqref{gammapoles} and let $(-\underline{\Delta}_{s})^{\sigma}$ be the fractional Laplacian defined in Theorem \ref{ThFracLap}. If
 $$
 w\in \bigcup_{\varepsilon>0}\mathcal{H}_{p}^{s+2+\frac{n+1}{p}+\varepsilon,\max\{\gamma+2,\frac{n+3}{2}\}+\varepsilon}(\mathbb{B})\oplus\mathbb{C}_{\omega} 
 $$ 
satisfies $w\geq \alpha>0$ on $\mathbb{B}$, for certain $\alpha>0$, then for each $\theta\in[0,\pi)$ there exists a $c>0$ such that $w(-\underline{\Delta}_{s})^{\sigma}+c\in \mathcal{R}(\theta)$.
\end{theorem}
\begin{proof}
By Definition \ref{MellSob} and \cite[Theorem 5.1]{Iz}, we have that the space $\mathcal{H}^{k,\widetilde{\gamma}}_{p}(\mathbb{B})$, $k\in\mathbb{N}_{0}$, $\widetilde{\gamma}\in\mathbb{R}$, admits an unconditional basis, see e.g. \cite[Section 4.1.b]{HNVW} for the notion of the unconditional basis of a Banach space. Therefore, by \cite[p. 527]{HP}, the space $\mathcal{H}^{k,\widetilde{\gamma}}_p(\mathbb{B})$ has the Pisier's {\em property $(\alpha)$}, see e.g. \cite[Definition 4.2.7]{PS} for the definition of the property $(\alpha)$. Then, by interpolation, i.e. by \cite[Theorem 3.11, Corollary 4.6]{KS} and \cite[Lemma 3.7]{RS2}, the space $\mathcal{H}^{\ell,\widetilde{\gamma}}_p(\mathbb{B})$, $\ell\in\mathbb{R}$, has the property $(\alpha)$ as well.

Denote $A=c_{0}-\underline{\Delta}_{s}$, $c_{0}>0$, and let $A^{\sigma}$ be defined by Theorem \ref{fracpower}. By \cite[Theorem 6.7]{SS1} we have that $A\in\mathcal{H}^{\infty}(\theta)$; in particular for each $\phi>0$ we have that $A\in \mathcal{BIP}(\phi)$. Hence, by \cite[Corollary 7.5 (b)]{KW0} the set $E=\{e^{-\phi|t|}A^{it}\, |\, t\in\mathbb{R}\}$ is $R$-bounded. Note that, in our situation, the $\gamma$-boundedness of the set $E$ (see \cite[Notations]{KW0} for the definition of this property), that is implied by \cite[Corollary 7.5 (b)]{KW0}, coincides with the $R$-boundedness, see \cite[Notations]{KW0} or the proof of \cite[Corollary 7.5 (b)]{KW0}. Moreover, by \cite[Lemma III.4.7.4]{Am} and \cite[Theorem 15.16]{KW1} we obtain that $A^{\sigma}\in \mathcal{BIP}(\sigma\phi)$ and $(A^{\sigma})^{it}=A^{i\sigma t}$, $t\in\mathbb{R}$. Therefore, we conclude that the set $\{e^{-\sigma\phi|t|}(A^{\sigma})^{it}\, |\, t\in\mathbb{R}\}$ is also $R$-bounded. Hence, \cite[Corollary 7.5 (a)]{KW0} (or \cite[Theorem 7.4]{Weis}) implies that $A^{\sigma}\in \mathcal{H}^{\infty}(\theta)$.

Let the operator $B:u\mapsto wu$, $u\in X_{0}^{s}$, which is bounded and invertible due to \cite[Lemma 3.3 and Lemma 6.2]{RS3}. Since in the situation of $B$ the path in formula \eqref{hgsta} can be chosen finite, by Kahane's contraction principle, see e.g. \cite[Proposition 2.5]{KW1}, we have that $B\in \mathcal{RH}^{\infty}(\theta)$, see also \cite[Theorem 4.5.4]{PS}. Moreover, $w\in \mathcal{H}_{p}^{\xi,\rho}(\mathbb{B})\oplus\mathbb{C}_{\omega}$ for some $\xi>s+2+\frac{n+1}{p}$ and $\rho>\max\{\gamma+2,\frac{n+3}{2}\}$, so that by \cite[Lemma 6.2]{RS3}
\begin{equation}\label{winverse}
(w+\mu)^{-1}\in \mathcal{H}_{p}^{\xi,\rho}(\mathbb{B})\oplus\mathbb{C}_{\omega} \quad \text{for each} \quad \mu\in S_{\theta}.
\end{equation}
In addition, from \eqref{sharpdom} $X_{\sigma}^{s}=\mathcal{H}_{p}^{s+2\sigma,\gamma+2\sigma}(\mathbb{B})\oplus\mathbb{C}_{\omega}$. Therefore, $B\mathcal{D}(A^{\sigma})\subseteq \mathcal{D}(A^{\sigma})$ and $(B+\mu)^{-1}\mathcal{D}(A^{\sigma})\subseteq \mathcal{D}(A^{\sigma})$, $\mu\in S_{\theta}$, due to \cite[Lemma 3.3]{RS3}. 

Choose $\nu\in(\max\{0,\sigma-1/2\},\sigma)$ such that $\gamma+2\nu-1\notin\cup_{j\in\mathbb{N}_{0}}\{ \pm\mu_{j}\}$. By \eqref{sharpdom} and \eqref{compare} we have that $\mathcal{D}(A^{\nu})=\mathcal{H}_{p}^{s+2\nu,\gamma+2\nu}(\mathbb{B})\oplus\mathbb{C}_{\omega}$, so that, if we denote by $\widetilde{B}$ the restriction of $B$ to $\mathcal{D}(A^{\nu})$, by \cite[Lemma 3.3]{RS3} we deduce that $\widetilde{B}\in \mathcal{L}(\mathcal{D}(A^{\nu}))$. Furthermore, $S_{\theta}\subset \rho(-\widetilde{B})$ and $(\widetilde{B}+\mu)^{-1}=(w+\mu)^{-1}$ when $\mu\in S_{\theta}$. Therefore, by \cite[Lemma 3.3]{RS3} we have
\begin{eqnarray}\nonumber
\lefteqn{\|(\widetilde{B}+\mu)^{-1}u\|_{\mathcal{D}(A^{\nu})}}\\\label{WKBD}
&\leq&C_{1}\|(\widetilde{B}+\mu)^{-1}\|_{\mathcal{H}_{p}^{\xi,\rho}(\mathbb{B})\oplus\mathbb{C}_{\omega}}\|u\|_{\mathcal{D}(A^{\nu})}=C_{1}|\mu|^{-1}\|(w\mu^{-1}+1)^{-1}\|_{\mathcal{H}_{p}^{\xi,\rho}(\mathbb{B})\oplus\mathbb{C}_{\omega}}\|u\|_{\mathcal{D}(A^{\nu})}
\end{eqnarray}
when $|\mu|\geq1$, for certain $C_{1}>0$. The set $\{w\mu^{-1}+1\, |\, \mu\in S_{\theta}, |\mu|\geq1\}$ is bounded in the space $\mathcal{H}_{p}^{\xi,\rho}(\mathbb{B})\oplus\mathbb{C}_{\omega}$, and moreover, there exists a $C_{2}>0$ such that $|w\mu^{-1}+1|>C_{2}$ when $\mu\in S_{\theta}$ and $|\mu|\geq1$. Hence, by \cite[Lemma 6.3]{RS3}, the set $\{\|(w\mu^{-1}+1)^{-1}\|_{\mathcal{H}_{p}^{\xi,\rho}(\mathbb{B})\oplus\mathbb{C}_{\omega}}\, |\, \mu\in S_{\theta}, |\mu|\geq1\}$ is also bounded and \eqref{WKBD} implies that $\widetilde{B}\in \mathcal{P}(\theta)$, i.e.
\begin{equation}\label{estnB}
\| (\widetilde{B}+\mu)^{-1}\|_{\mathcal{L}(\mathcal{D}(A^{\nu}))}\leq \frac{C_{3}}{1+|\mu|}, \quad \mu\in S_{\theta},
\end{equation}
for certain $C_{3}>0$; when $\mu\in S_{\theta}$, $|\mu|\leq1$ in \eqref{estnB}, the norm $\|(\widetilde{B}+\mu)^{-1}\|_{\mathcal{L}(\mathcal{D}(A^{\nu}))}$ is estimated by $\|(w+\mu)^{-1}\|_{\mathcal{H}_{p}^{\xi,\rho}(\mathbb{B})\oplus\mathbb{C}_{\omega}}$, which is bounded due to \cite[Lemma 6.3]{RS3}.

By \cite[Theorem III.4.6.13]{Am} we have $A^{\nu}=(A^{\sigma})^{\nu/\sigma}$ so that, from Lemma \ref{lemmadecay}, we infer
$$
\|A^{\nu}(A^{\sigma}+\lambda)^{-1}\|_{\mathcal{L}(X_{0})}\leq\frac{C_{4}}{1+|\lambda|^{1-\frac{\nu}{\sigma}}}, \quad \lambda\in S_{\theta},
$$
for certain $C_{4}>0$. Furthermore, by Lemma \ref{commlemma} and the regularity of $w$, we have that $[A^{\sigma},B]A^{-\nu}\in\mathcal{L}(X_{0}^{s})$. Hence, by taking into account \eqref{estnB}, we estimate
\begin{eqnarray*}
\lefteqn{\|[A^{\sigma},(B+\mu)^{-1}](A^{\sigma}+\lambda)^{-1}\|_{\mathcal{L}(X_{0}^{s})}}\\
&=&\|(B+\mu)^{-1}[A^{\sigma},B](B+\mu)^{-1}(A^{\sigma}+\lambda)^{-1}\|_{\mathcal{L}(X_{0}^{s})}\\
&=&\|(B+\mu)^{-1}[A^{\sigma},B]A^{-\nu}A^{\nu}(\widetilde{B}+\mu)^{-1}A^{-\nu}A^{\nu}(A^{\sigma}+\lambda)^{-1}\|_{\mathcal{L}(X_{0}^{s})}\\
&\leq&\|(B+\mu)^{-1}\|_{\mathcal{L}(X_{0}^{s})}\|[A^{\sigma},B]A^{-\nu}\|_{\mathcal{L}(X_{0}^{s})}\|A^{\nu}(\widetilde{B}+\mu)^{-1}A^{-\nu}\|_{\mathcal{L}(X_{0}^{s})}\|A^{\nu}(A^{\sigma}+\lambda)^{-1}\|_{\mathcal{L}(X_{0}^{s})}\\
&\leq&\frac{C_{5}}{(1+|\mu|^{1+\eta})(1+|\lambda|^{1-\frac{\nu}{\sigma}})},
\end{eqnarray*}
for all $\lambda,\mu\in S_{\theta}$, all $\eta\in(0,1)$ and certain $C_{5}>0$. We conclude that $A^{\sigma}$ and $B$ satisfy the Da Prato and Grisvard commutation condition \eqref{DaPratoGrisv} and, by Theorem \ref{HaHiTh}, there exists a $c>0$ such that $A^{\sigma}B+c\in \mathcal{H}^{\infty}(\theta)$, where we have used the fact that $\theta$ can be chosen arbitrary close to $\pi$. In particular, see \cite[Theorem 4.4.5]{PS}, $A^{\sigma}B+c\in \mathcal{R}(\theta)$. 

Recall that due to \eqref{winverse}, $w^{-1}\in\mathcal{H}_{p}^{\xi,\rho}(\mathbb{B})\oplus\mathbb{C}_{\omega}$, so that by \cite[Lemma 3.3]{RS3} multiplication by $w$ or $w^{-1}$ induces a bounded map on $X_{0}^{s}$. Moreover, by the regularity of $w$ and \cite[Lemma 3.3]{RS3}, we have that $BA^{\sigma}\in \mathcal{L}(X_{\sigma}^{s},X_{0}^{s})$. Hence, by the formula $BA^{\sigma}+c+\lambda=B(A^{\sigma}B+c+\lambda)B^{-1}$, $\lambda\in S_{\theta}$, we deduce that $S_{\theta}\subset \rho(-(BA^{\sigma}+c))$ and $(BA^{\sigma}+c+\lambda)^{-1}=B(A^{\sigma}B+c+\lambda)^{-1}B^{-1}$ when $\lambda\in S_{\theta}$. This resolvent representation together with \cite[Lemma 3.3]{RS3} and the definition of $R$-sectoriality imply that $BA^{\sigma}+c\in \mathcal{R}(\theta)$. The result then follows by Step 2 of the proof of Theorem \ref{wperturb}, i.e. the same argument is applicable to the case of $s>0$.
\end{proof}

Before we proceed to the proof of the main fractional porous medium equation result, we recall certain embedding properties of the real interpolation between Mellin-Sobolev spaces.

\begin{corollary}\label{corembofint}
Let $p,q\in(1,\infty)$, $s\geq0$, $\gamma$ be as in \eqref{gammaweight} and $\sigma\in(0,1)$. The following embeddings hold
\begin{eqnarray}\nonumber
\lefteqn{\bigcup_{\varepsilon>0} \mathcal{H}_{p}^{s+2\sigma-\frac{2\sigma}{q}+\varepsilon,\gamma+2\sigma-\frac{2\sigma}{q}+\varepsilon}(\mathbb{B})\oplus\mathbb{C}_{\omega}}\\\label{embedtoH}
&&\hookrightarrow (X_{\sigma}^{s},X_{0}^{s})_{\frac{1}{q},q}\hookrightarrow \bigcap_{\varepsilon>0} \mathcal{H}_{p}^{s+2\sigma-\frac{2\sigma}{q}-\varepsilon,\gamma+2\sigma-\frac{2\sigma}{q}-\varepsilon}(\mathbb{B})\oplus\mathbb{C}_{\omega}. 
\end{eqnarray}
If in addition
$$
s+2\sigma-\frac{2\sigma}{q}>\frac{n+1}{p} \quad \text{and} \quad \gamma+2\sigma-\frac{2\sigma}{q}>\frac{n+2}{2}, 
$$
then 
\begin{equation}\label{embedtoC}
 \bigcap_{\varepsilon>0} \mathcal{H}_{p}^{s+2\sigma-\frac{2\sigma}{q}-\varepsilon,\gamma+2\sigma-\frac{2\sigma}{q}-\varepsilon}(\mathbb{B})\oplus\mathbb{C}_{\omega}\hookrightarrow C(\mathbb{B}).
\end{equation}
\end{corollary}
\begin{proof}
By reiteration, see e.g. \cite[(I.2.5.2), (I.2.8.1), (I.2.8.2) and (I.2.9.6)]{Am}, we have $(X_{\sigma}^{s},X_{0}^{s})_{1/q,q}=(X_{0}^{s},X_{1}^{s})_{\sigma(1-1/q),q}$, so that the first embedding follows by \cite[(I.2.5.4)]{Am} and \cite[Lemma 5.2]{RS3}. The second embedding follows by \cite[Corollary 2.9]{RS2}.
\end{proof}

\subsubsection*{Proof of Theorem 1.3}

As a first step we apply Theorem \ref{ClementLi} to 
\begin{eqnarray}\label{FPME3}
w'(t)+mw^{\frac{m-1}{m}}(-\Delta)^{\sigma}w(t) &=& 0, \quad t\in(0,T),\\\label{FPME4}
w(0) &=& w_{0}=u_{0}^{m},
\end{eqnarray}
with $A(\cdot)=m(\cdot)^{\frac{m-1}{m}}(-\underline{\Delta}_{s})^{\sigma}$, $s\geq0$, and the Banach couple $X_{0}^{s}$, $X_{\sigma}^{s}$. If $u_{0}$ is as in \eqref{u01} or \eqref{u02}, then by \cite[Lemma 6.2]{RS3} and \eqref{embedtoH} we have respectively that $w_{0}\in \mathcal{H}_{p}^{2\sigma-\frac{2\sigma}{q}-\varepsilon,\gamma+2\sigma-\frac{2\sigma}{q}-\varepsilon}(\mathbb{B})\oplus\mathbb{C}_{\omega}$ for all $\varepsilon>0$ small enough or $w_{0}\in \mathcal{H}_{p}^{\nu+2+\frac{n+1}{p}+\widetilde{\varepsilon},\max\{\gamma+2,\frac{n+3}{2}\}+\widetilde{\varepsilon}}(\mathbb{B})\oplus\mathbb{C}_{\omega}$ for some $\widetilde{\varepsilon}>0$. Therefore, the maximal $L^q$-regularity of the linearized term follows by Theorem \ref{KaWeTh}, Theorem \ref{wperturb}, Theorem \ref{wperturbsect} and \eqref{embedtoH}. 

The Lipschitz continuity of $A(\cdot)$ follows similarly to \cite[(6.20)]{RS3}. More precisely, let $B_{r}$ be an open ball in $(X_{\sigma}^{s},X_{0}^{s})_{1/q,q}$ of radius $r>0$ centered at $w_{0}$. Due to \eqref{embedtoC}, choose $r>0$ sufficiently small and let $\Gamma$ be a finite closed path in $\{\lambda\in\mathbb{C}\, |\, \mathrm{Re}(\lambda)<0\}$ that surrounds $\cup_{v\in B_{r}}\mathrm{Ran}(-v)$. For each $\eta\in\mathbb{R}$ we have
\begin{equation}\label{fraccontins}
w_{1}^{\eta}-w_{2}^{\eta}=(w_{2}-w_{1})\frac{1}{2\pi i}\int_{\Gamma}(-\lambda)^{\eta}(w_{1}+\lambda)^{-1}(w_{2}+\lambda)^{-1}d\lambda, \quad w_{1},w_{2}\in B_{r}.
\end{equation}
Therefore, by the above formula, \cite[Corollary 3.2, Corollary 3.3, Lemma 6.2, Lemma 6.3]{RS3} and \eqref{embedtoH} we estimate
\begin{eqnarray}\nonumber
\lefteqn{\|w_{1}^{\eta}(-\underline{\Delta}_{s})^{\sigma}-w_{2}^{\eta}(-\underline{\Delta}_{s})^{\sigma}\|_{\mathcal{L}(X_{\sigma}^{s},X_{0}^{s})}\leq C_{1}\|w_{1}^{\eta}-w_{2}^{\eta}\|_{\mathcal{L}(X_{0}^{s})}}\\\nonumber
&\leq&C_{2}\|w_{1}-w_{2}\|_{\mathcal{H}_{p}^{s+2\sigma-\frac{2\sigma}{q}-\varepsilon,\gamma+2\sigma-\frac{2\sigma}{q}-\varepsilon}(\mathbb{B})\oplus\mathbb{C}_{\omega}}\\\nonumber
&&\times\int_{\Gamma}|\lambda|^{\eta}\Big(\|(w_{1}+\lambda)^{-1}\|_{\mathcal{H}_{p}^{s+2\sigma-\frac{2\sigma}{q}-\varepsilon,\gamma+2\sigma-\frac{2\sigma}{q}-\varepsilon}(\mathbb{B})\oplus\mathbb{C}_{\omega}}\\\label{Lipschest}
&&\times\|(w_{2}+\lambda)^{-1}\|_{\mathcal{H}_{p}^{s+2\sigma-\frac{2\sigma}{q}-\varepsilon,\gamma+2\sigma-\frac{2\sigma}{q}-\varepsilon}(\mathbb{B})\oplus\mathbb{C}_{\omega}}\Big)d\lambda\leq C_{3}\|w_{1}-w_{2}\|_{(X_{\sigma}^{s},X_{0}^{s})_{\frac{1}{q},q}},
\end{eqnarray}
for certain $C_{1}, C_{2}, C_{3}>0$ and all $\varepsilon>0$ sufficiently small. By choosing $\eta=\frac{m-1}{m}$ we conclude that there exists a $T>0$ and a unique 
\begin{eqnarray}\nonumber
\lefteqn{w\in W^{1,q}(0,T;X_{0}^{s})\cap L^{q}(0,T; X_{\sigma}^{s})}\\\label{reqw}
&&\hookrightarrow \bigcap_{\varepsilon>0}C([0,T]; \mathcal{H}_{p}^{s+2\sigma-\frac{2\sigma}{q}-\varepsilon,\gamma+2\sigma-\frac{2\sigma}{q}-\varepsilon}(\mathbb{B})\oplus\mathbb{C}_{\omega})\hookrightarrow C([0,T];C(\mathbb{B}))
\end{eqnarray}
solving \eqref{FPME3}-\eqref{FPME4}, where $s=0$ refers to \eqref{u01} and $s=\nu$ to \eqref{u02}; note that, by uniqueness, the solution for $s=\nu$ coincides (possibly in a smaller interval $[0,T]$) with the solution for $s=0$. Moreover, in \eqref{reqw} we have used \eqref{interpemb} and Corollary \ref{corembofint}.

Similarly to \eqref{Lipschest}, by Cauchy's integral formula, for $\varepsilon>0$ sufficiently small we have
\begin{eqnarray}\nonumber
\lefteqn{\|w_{1}^{\eta}(-\underline{\Delta}_{s})^{\sigma}-w_{2}^{\eta}(-\underline{\Delta}_{s})^{\sigma}-\eta(w_{1}-w_{2})w_{1}^{\eta-1}(-\underline{\Delta}_{s})^{\sigma}\|_{\mathcal{L}(X_{\sigma}^{s},X_{0}^{s})}}\\\nonumber
&\leq&C_{4}\|w_{1}^{\eta}-w_{2}^{\eta}-\eta(w_{1}-w_{2})w_{1}^{\eta-1}\|_{\mathcal{H}_{p}^{s+2\sigma-\frac{2\sigma}{q}-\varepsilon,\gamma+2\sigma-\frac{2\sigma}{q}-\varepsilon}(\mathbb{B})\oplus\mathbb{C}_{\omega}}\\\nonumber
&\leq&\frac{C_{4}}{2\pi}\int_{\Gamma}|\lambda|^{\eta}\|(w_{1}-w_{2})^{2}(w_{1}+\lambda)^{-2}(w_{2}+\lambda)^{-1}\|_{\mathcal{H}_{p}^{s+2\sigma-\frac{2\sigma}{q}-\varepsilon,\gamma+2\sigma-\frac{2\sigma}{q}-\varepsilon}(\mathbb{B})\oplus\mathbb{C}_{\omega}}d\lambda\\\label{freschetder}
&\leq&C_{5} \|w_{1}-w_{2}\|_{(X_{\sigma}^{s},X_{0}^{s})_{\frac{1}{q},q}}^{2} \int_{\Gamma}|\lambda|^{\eta}\|(w_{1}+\lambda)^{-1}\|_{(X_{\sigma}^{s},X_{0}^{s})_{\frac{1}{q},q}}^{2}\|(w_{2}+\lambda)^{-1}\|_{(X_{\sigma}^{s},X_{0}^{s})_{\frac{1}{q},q}}d\lambda,
\end{eqnarray}
for certain $C_{4}, C_{5}>0$. Since $\eta$ was arbitrary, by applying successively \eqref{freschetder}, we obtain in particular that $A(\cdot)\in C^{\infty}(B_{r};\mathcal{L}(X_{\sigma}^{0},X_{0}^{0}))$. 

Moreover, if we choose $r_{0}>0$ sufficiently small and restrict $r\in(0,r_{0})$, then the constant bound $C_{3}$ in \eqref{Lipschest} can be chosen independently of $r\in(0,r_{0})$. Therefore, by taking $r>0$ small enough, due to \eqref{Lipschest} and the perturbation result \cite[Proposition 4.4.2]{PS}, for each $v\in B_{r}$ we have that $A(v)$ has maximal $L^{q}$-regularity. Thus, by \eqref{interpemb} we can restrict to a sufficiently small $T>0$ such that $w(t)\in B_{r}$, $t\in[0,T]$, and in particular each $A(w(t))$, $t\in[0,T]$, has maximal $L^{q}$-regularity. Then by \cite[Theorem 5.2.1]{PS}, in addition to \eqref{reqw} we have 
\begin{equation}\label{wsmooth}
w\in C^{\infty}((0,T);X_{\sigma}^{0}).
\end{equation}

By letting $u=w^{1/m}$, from \eqref{FPME3}-\eqref{FPME4} we see that $u$ satisfies the original equation \eqref{FPME1}-\eqref{FPME2}. Hence, it suffices to show that $u$ also satisfies the regularity \eqref{reqw} and \eqref{wsmooth}. By \eqref{fraccontins}, \cite[Corollary 3.2, Lemma 6.2, Lemma 6.3]{RS3} and \eqref{reqw}, similarly to the estimate \eqref{Lipschest}, we have that 
\begin{equation}\label{firsttegu}
u, u^{1-m}, u^{\frac{1-m}{m}}\in \bigcap_{\varepsilon>0}C([0,T]; \mathcal{H}_{p}^{s+2\sigma-\frac{2\sigma}{q}-\varepsilon,\gamma+2\sigma-\frac{2\sigma}{q}-\varepsilon}(\mathbb{B})\oplus\mathbb{C}_{\omega}).
\end{equation}
Hence, $u\in L^{q}(0,T; X_{0}^{s})$. In addition, by the formula $\partial_{t}u=m^{-1}w^{\frac{1-m}{m}}\partial_{t}w$, \eqref{reqw}, \eqref{firsttegu} and \cite[Corollary 3.3]{RS3} we also have $u'\in L^{q}(0,T; X_{0}^{s})$, so that $u\in W^{1,q}(0,T;X_{0}^{s})$. 

Recall that from \eqref{sharpdom} and \cite[Lemma 6.2]{RS3}, we have $u(t)\in X_{\sigma}^{s}$ for almost all $t\in[0,T]$. Therefore by \eqref{FPME1}, \eqref{firsttegu} and \cite[Corollary 3.2]{RS3}, for $\varepsilon>0$ small enough we estimate 
\begin{eqnarray*}
\lefteqn{\int_{0}^{T}\|u(t)\|_{X_{\sigma}^{s}}^{q}dt \leq C_{6} \int_{0}^{T}\|u^{1-m}(t)\|_{\mathcal{H}_{p}^{s+2\sigma-\frac{2\sigma}{q}-\varepsilon,\gamma+2\sigma-\frac{2\sigma}{q}-\varepsilon}(\mathbb{B})\oplus\mathbb{C}_{\omega}}^{q}\|u^{m}(t)\|_{X_{\sigma}^{s}}^{q}dt}\\
&\leq&C_{7}\int_{0}^{T}(\|(-\underline{\Delta}_{s})^{\sigma}u^{m}(t)\|_{X_{0}^{s}}+\|u(t)\|_{X_{0}^{s}})^{q}dt=C_{7}\int_{0}^{T}(\|u'(t)\|_{X_{0}^{s}}+\|u(t)\|_{X_{0}^{s}})^{q}dt,
\end{eqnarray*}
for certain $C_{6}, C_{7}>0$, and by the triangle inequality we conclude that $u\in L^{q}(0,T; X_{\sigma}^{s})$. Finally, in the case of $s=0$, \eqref{fraccontins}, \eqref{wsmooth}, \cite[Lemma 3.2, Lemma 6.2 and Lemma 6.3]{RS3} imply that $u\in C^{\infty}((0,T);X_{\sigma}^{0})$.
\mbox{\ } \hfill $\square$

\section{Appendix}

In this section we collect some elementary abstract results we have used previously. We recall first the following decay property of the resolvent of a sectorial operator.

\begin{lemma}\label{lemmadecay}
If $A\in\mathcal{P}(\theta)$, $\theta\in[0,\pi)$, in $X_{0}$, then for any $\sigma\in[0,1]$ there exists a $C>0$, depending only on $\theta$, the sectorial bound of $A$ and $\sigma$, such that 
$$
\|A^{\sigma}(A+\lambda)^{-1}\|_{\mathcal{L}(X_{0})}\leq\frac{C}{1+|\lambda|^{1-\sigma}}, \quad \lambda\in S_{\theta}.
$$
\end{lemma}
\begin{proof}
It is sufficient to consider the case of $\sigma\in(0,1)$ and show the estimate for $|\lambda|>1$. Recall that $A\in\mathcal{P}(\phi)$, for some $\phi\in(\theta,\pi)$. Thus, for any $\lambda\in S_{\theta}$, $|\lambda|>1$, by \eqref{cp} and Cauchy's theorem we have
\begin{eqnarray*}
\lefteqn{A^{\sigma}(A+\lambda)^{-1}=\frac{1}{2\pi i}A\int_{\Gamma_{\phi}}(-z)^{\sigma-1}(A+z)^{-1}(A+\lambda)^{-1}dz}\\
&=&\frac{1}{2\pi i}A\int_{\Gamma_{\phi}}\frac{(-z)^{\sigma-1}}{\lambda-z}((A+z)^{-1}-(A+\lambda)^{-1})dz\\
&=&\frac{1}{2\pi i}A\int_{\Gamma_{\phi}}\frac{(-z)^{\sigma-1}}{\lambda-z}(A+z)^{-1}dz-\frac{1}{2\pi i}A(A+\lambda)^{-1}\int_{\Gamma_{\phi}}\frac{(-z)^{\sigma-1}}{\lambda-z}dz\\
&=&\frac{1}{2\pi i}\int_{\Gamma_{\phi}}\frac{(-z)^{\sigma-1}}{\lambda-z}(A+z-z)(A+z)^{-1}dz\\
&=&\frac{1}{2\pi i}\int_{\Gamma_{\phi}}\frac{z(-z)^{\sigma-1}}{z-\lambda}(A+z)^{-1}dz.
\end{eqnarray*}
Hence,
\begin{eqnarray*}
\lambda^{1-\sigma}A^{\sigma}(A+\lambda)^{-1}&=&\frac{1}{2\pi i}\int_{-\Gamma_{\phi}}\frac{(\frac{\lambda}{z})^{1-\sigma}}{1+\frac{\lambda}{z}}(A-z)^{-1}dz\\
&=&\frac{1}{2\pi i}\int_{-\Gamma_{\phi}}\frac{(\frac{\lambda}{|\lambda|})^{1-\sigma}w^{-\sigma}}{1+\frac{\lambda}{|\lambda|}w}(\frac{|\lambda|}{w}-A)^{-1}\frac{|\lambda|}{w}dw,
\end{eqnarray*}
and the estimate follows.
\end{proof}

We end up with a mapping property of the coefficients of the resolvent's Laurent expansion.

\begin{lemma}\label{expandatzero}
Let $A:\mathcal{D}(A)\rightarrow X_{0}$ be a closed linear operator in $X_{0}$ such that $0\notin \rho(-A)$. Assume that there exists some neighbourhood $U$ of $0$ such that 
$$
(A+\lambda)^{-1}=\lambda^{-1}B_{-1}+\sum_{k=0}^{\infty}\lambda^{k} B_{k} , \quad \lambda\in U\backslash\{0\},
$$
for some $B_{k}\in\mathcal{L}(X_{0})$, $k\in\{-1\}\cup\mathbb{N}_{0}$, where the series converges absolutely in the $\mathcal{L}(X_{0})$-topology. Then, $B_{k}\in \mathcal{L}(X_{0},\mathcal{D}(A))$, $k\in\{-1\}\cup\mathbb{N}_{0}$, and
\begin{equation}\label{edtar}
AB_{-1}=0, \quad AB_{0}=I-B_{-1}, \quad AB_{k}=-B_{k-1}, \quad k\geq1.
\end{equation}
Consequently, the series converges absolutely in the $\mathcal{L}(X_{0},\mathcal{D}(A))$-topology as well.
\end{lemma}
\begin{proof}
If $u\in X_{0}$, then 
\begin{equation}\label{hagtr}
\lambda(A+\lambda)^{-1}u \rightarrow B_{-1}u \quad \text{as} \quad \lambda\rightarrow0
\end{equation}
and 
$$
\lambda A(A+\lambda)^{-1}u=\lambda(I-\lambda(A+\lambda)^{-1})u\rightarrow 0 \quad \text{as} \quad \lambda\rightarrow0.
$$
Therefore, by the closedness of $A$ we conclude that 
\begin{equation}\label{step1}
B_{-1}u\in \mathcal{D}(A) \quad \text{and}\quad AB_{-1}u=0.
\end{equation}

Similarly, we have
$$
\lambda^{-1}(\lambda(A+\lambda)^{-1}u-B_{-1}u)\rightarrow B_{0}u \quad \text{as} \quad \lambda\rightarrow0.
$$
Also, by \eqref{hagtr} and \eqref{step1} we obtain 
$$
\lambda^{-1}A(\lambda(A+\lambda)^{-1}u-B_{-1}u)=u-\lambda(A+\lambda)^{-1}u\rightarrow u-B_{-1}u\quad \text{as} \quad \lambda\rightarrow0.
$$
Hence, the closedness of $A$ implies that 
\begin{equation}\label{step2}
B_{0}u\in \mathcal{D}(A) \quad \text{and}\quad AB_{0}u=u-B_{-1}u.
\end{equation}

Moreover,
$$
\lambda^{-1}(\lambda^{-1}(\lambda(A+\lambda)^{-1}u-B_{-1}u)-B_{0}u)\rightarrow B_{1}u\quad \text{as} \quad \lambda\rightarrow0,
$$
and by \eqref{step1}-\eqref{step2}
\begin{eqnarray*}
\lefteqn{\lambda^{-1}A(\lambda^{-1}(\lambda(A+\lambda)^{-1}u-B_{-1}u)-B_{0}u)}\\
&=&\lambda^{-1}((u-\lambda(A+\lambda)^{-1}u)-u+B_{-1}u)=-(A+\lambda)^{-1}u+\lambda^{-1}B_{-1}u\rightarrow -B_{0}u\quad \text{as} \quad \lambda\rightarrow0.
\end{eqnarray*}
Thus, $B_{1}u\in \mathcal{D}(A)$ and $AB_{1}u=-B_{0}u$.

Assume that for some $\nu\in\mathbb{N}_{0}$, $\nu\geq1$, we have $B_{k}\in \mathcal{L}(X_{0},\mathcal{D}(A))$, $k\leq \nu$, and in addition \eqref{edtar} holds for all $k\leq \nu$. We have
$$
\lambda^{-1}(\lambda^{-1}(\cdots\lambda^{-1}(\lambda^{-1}(\lambda(A+\lambda)^{-1}-B_{-1})-B_{0})\cdots)-B_{\nu})u\rightarrow B_{\nu+1}u\quad \text{as} \quad \lambda\rightarrow0.
$$
Furthermore, by \eqref{edtar} we obtain
\begin{eqnarray*}
\lefteqn{A\lambda^{-1}(\lambda^{-1}(\cdots\lambda^{-1}(\lambda^{-1}(\lambda(A+\lambda)^{-1}-B_{-1})-B_{0})\cdots)-B_{\nu})u}\\
&=&\lambda^{-\nu-1}A(\lambda^{-1}(\lambda(A+\lambda)^{-1}-B_{-1})-B_{0})u-\sum_{k=1}^{\nu}\lambda^{k-\nu-1}AB_{k}u\\
&=&\lambda^{-\nu-1}B_{-1}u-\lambda^{-\nu}(A+\lambda)^{-1}u +\sum_{k=0}^{\nu-1}\lambda^{k-\nu}B_{k}u\rightarrow -B_{\nu}u \quad \text{as} \quad \lambda\rightarrow0.
\end{eqnarray*}
Hence, by the closedness of $A$, $B_{\nu+1}u\in \mathcal{D}(A)$ and $AB_{\nu+1}u=-B_{\nu}u$.
\end{proof}

\begin{remark}\label{generlemmaresv}
By following the same proof, Lemma \ref{expandatzero} can be generalized from $0$ to an arbitrary pole $\lambda_{0}$ of order $\mu\in\mathbb{N}$. More precisely, let $A:\mathcal{D}(A)\rightarrow X_{0}$ be a closed linear operator in $X_{0}$ such that $\lambda_{0}\notin \rho(-A)$. Assume that there exists some neighbourhood $U$ of $\lambda_{0}$ such that 
$$
(A+\lambda)^{-1}=\sum_{k=-\mu}^{\infty}(\lambda-\lambda_{0})^{k} B_{k} , \quad \lambda\in U\backslash\{\lambda_{0}\},
$$
for some $\mu\in\mathbb{N}$ and certain $B_{k}\in\mathcal{L}(X_{0})$, $k\geq-\mu$, where the series converges absolutely in the $\mathcal{L}(X_{0})$-topology. Then, $B_{k}\in \mathcal{L}(X_{0},\mathcal{D}(A))$, $k\geq-\mu$, and 
$$
AB_{-\mu}+\lambda_{0}B_{-\mu}=0,\quad AB_{0}+B_{-1}+\lambda_{0} B_{0}=I, \quad AB_{k}+B_{k-1}+\lambda_{0} B_{k}=0, \quad k\geq1-\mu, \, k\neq 0.
$$
Consequently, the series converges absolutely in the $\mathcal{L}(X_{0},\mathcal{D}(A))$-topology as well.
\end{remark}

\begin{remark}
Lemma \ref{expandatzero} and Remark \ref{generlemmaresv} can be alternatively shown by using first Cauchy's integral formula and then the identity $A(A+\lambda)^{-1}=I-\lambda(A+\lambda)^{-1}$, $\lambda\in\rho(-A)$.
\end{remark}

\end{document}